\documentclass[12pt,a4paper]{amsart}
\usepackage{geometry}
\usepackage{latexsym}
\usepackage{bm}
\usepackage{amsmath}
\usepackage{amssymb}
\usepackage{amsmath,amscd}
\usepackage{tabls}
\usepackage{url}
\usepackage[utf8]{inputenc}
\usepackage{amsfonts}
\usepackage{amssymb}
\usepackage{mathrsfs}
\usepackage{amsmath}
\usepackage{amsthm}
\usepackage{amscd}
\usepackage{fancyhdr}
\usepackage{enumerate}
\usepackage{paralist}
\usepackage{xypic}
\usepackage{graphicx}
\usepackage{cite}
\usepackage{lipsum}
\usepackage{footmisc}
\usepackage[all,cmtip]{xy}

\theoremstyle{plain}

\newtheorem{Thm}[equation]{Theorem}
\newtheorem{lem}[equation]{Lemma}
\newtheorem{prop}[equation]{Proposition}
\newtheorem{rem}[equation]{Remark}

\newtheorem{conj}[equation]{Conjecture}

\numberwithin{equation}{section}

\title{Motivic multiple zeta values relative to $\mu_2$}
\author{Zhongyu Jin}
\email{zyjin@pku.edu.cn}
\address{Zhongyu Jin \\School of Mathematical Sciences,
        Peking University,
         Beijing, China}

\author{Jiangtao Li}
\email{ljt-math@pku.edu.cn}
\address{Jiangtao Li \\School of Mathematical Sciences,
        Peking University,
         Beijing, China}

\begin{document}
\maketitle
\begin{abstract}
We establish a short exact sequence about depth-graded motivic double zeta values of even weight relative to $\mu_2$. We find a basis for the depth-graded motivic double zeta values relative to $\mu_2$ of even weight and a basis for the depth-graded motivic triple zeta values relative to $\mu_2$ of odd weight. As an application of our main results, we prove Kaneko and Tasaka's conjectures about the  sum odd double zeta values and the  classical double zeta values.  We also prove an analogue of Kaneko and Tasaka's conjecture in depth three. At last we formulate a conjecture which is related to  sum odd multiple zeta values in higher depth.
\end{abstract}
\let\thefootnote\relax\footnotetext{
2010 $\mathnormal{Mathematics} \;\mathnormal{Subject}\;\mathnormal{Classification}$. Primary 11F32, Secondary 11F67.\\
$\mathnormal{Keywords:}$  Multiple zeta values, Period polynomial, Mixed Tate motives. }

\maketitle

\section{Introduction}

Multiple zeta values are defined by the following convergent series:
$$\zeta (n_{1},\cdots,n_{r})=\sum_{0<k_{1}<\cdots<k_{r}} \frac{1}{k_{1}^{n_{1}}\cdots k_{r}^{n_{r}}},\ n_{1},...,n_{r-1}>0,\ n_{r}>1.$$
In particularly, when $r=1$ they are  the classical Riemann zeta values. We call $r$  the depth, and $N=n_{1}+\cdots+n_{r}$  the weight of the above multiple zeta value.

Denote by $\mathcal{Z}_N$ the $\mathbb{Q}$-vector space generated by all the weight $N$ multiple zeta values for $N> 0$, and $\mathcal{Z}_0=\mathbb{Q}$. Then
\[\mathcal{Z}=\bigoplus_{N\geq 0}\mathcal{Z}_N\]
  is  a graded algebra with the shuffle product. There is a depth filtration $\mathcal{D}$ on $\mathcal{Z}_N$
$$ \mathcal{D}_{r}\mathcal{Z}_N=\langle\zeta (n_{1},\cdots,n_{k})\in \mathcal{Z}_N, k\leq r\rangle_{\mathbb{Q}},$$
where for $x_i\in \mathbb{R},i\in A$, $\langle x_i, i\in A\rangle_{\mathbb{Q}}$ means the $\mathbb{Q}$-linear subspace generated by $x_i,i\in A$ in $\mathbb{R}$.

The double zeta values are the subspace $\mathcal{DZ}$ of $\mathcal{Z}$ spanned by elements which have depth $2$.  Gangl, Kaneko, Zagier \cite{kan} found an interesting connection between period polynomials of $SL_2(\mathbb{Z})$ and the double shuffle relations among $\mathcal{DZ}$.

Brown  \cite{ref2}   defined the motivic multiple zeta values algebra $\mathcal{H}$, its elements can be written as  $\mathbb{Q}$-linear combinations of motivic multiple zeta values $\zeta^{\mathfrak{m}} (n_{1},...,n_{r})$. Motivic multiple zeta values also satisfy the double shuffle relations by the work of Soud\`{e}res \cite{ref10}. There is a surjective graded algebra homomorphism:
$$\eta: \mathcal{H} \rightarrow \mathcal{Z},$$
$$\zeta^{\mathfrak{m}} (n_{1},\cdots,n_{r}) \mapsto \zeta (n_{1},\cdots,n_{r}).$$
Brown proved that   the set $\{\zeta^{\mathfrak{m}} (n_{1},\cdots,n_{r}), n_{i}\in {2,3}\}$ is a basis for non-zero weight subspace of  $\mathcal{H}$, thus proved that every multiple zeta value is a $\mathbb{Q}$-linear combination of ${\zeta (n_{1},\cdots,n_{r}),\ n_{i}\in {2,3}}$ (Hoffman's conjecture C in \cite{hof}).
It suggests that we can study the multiple zeta values by studying these motivic multiple zeta values.

Since  motivic multiple zeta values also satisfy the double shuffle relations. By Gangl, Kaneko, Zagier's results, there are  also period polynomial relations among motivic double zeta values of even weight.   This fact was reintepreted as a short exact sequence which inlvoves motivic double zeta values of even weight with a slight modification and period polynomials in the second author's paper \cite{ref8}.

 Futhermore, the second author  \cite{ref8} proposed two exact sequence conjectures which relate the depth-graded version of motivic multiple zeta values and period polynomials of $SL_2(\mathbb{Z})$. The second author verified the two conjectures in low depth. Besides, the second author and Liu \cite{ref9}  established a short exact sequence about motivic double zeta values of odd weight.

After Brown, Glanois \cite{ref5} defined the motivic multiple zeta values relative to $\mu_{N}$, which is a generalization of $\mathcal{H}$ in the cyclotomic field , where $\{\mu_{N}\} $ is the set of all $N^{th}$ roots of unity.  She gave a basis of motivic multiple zeta values relative to $\mu_N$ for $N=2,3,4,6,8$.
We will give a brief introduction to  Glanois' work in the next section in the case of $N=2$.

Ma \cite{ma} studied motivic double zeta values relative to $\mu_N$ for $N=2,3$. He found various connections between some special matrices which come from motivic Galois action on motivic double zeta values relative to $\mu_N$, Hecke operators and newforms of $\Gamma_0(N)$ for $N=2,3$.

In the rest of this paper, we only consider the motivic multiple zeta values relative to $\mu_2$, and  denote by $\mathcal{H}$ the $\mathbb{Q}$-algebra generated by them rather than the motivic multiple zeta  algebra of Brown for convenience.

For positive integers $n_{1}\geq 1, n_{2}\geq 2$,  define
$$\zeta^{o} (n_{1},n_{2})=\sum_{0<k_{1}<k_{2}, \mathrm{odd}} \frac{1}{k_{1}^{n_{1}}k_{2}^{n_{2}}}.$$
It is obvious to see
$$\zeta^{o} (n_{1},n_{2})=\frac{1}{4}\left(\zeta (n_{1},n_{2})-\zeta (\overline{n_{1}},n_{2})-\zeta (n_{1},\overline{n_{2}})+\zeta (\overline{n_{1}},\overline{n_{2}})\right),$$
where
$$\zeta (\overline{n_{1}},n_{2})=\sum_{0<k_{1}<k_{2}}\frac{ (-1)^{k_{1}}}{k_{1}^{n_{1}}k_{2}^{n_{2}}},$$
$$\zeta (n_{1},\overline{n_{2}})=\sum_{0<k_{1}<k_{2}} \frac{(-1)^{k_{2}}}{k_{1}^{n_{1}}k_{2}^{n_{2}}},$$
$$\zeta (\overline{n_{1}},\overline{n_{2}})=\sum_{0<k_{1}<k_{2}} \frac{(-1)^{k_{1}+k_{2}}}{k_{1}^{n_{1}}k_{2}^{n_{2}}}.$$

Denote by $\mathcal{Z}^{2}$ the $\mathbb{Q}$-vector space generated by
$$\zeta \binom{n_1,\cdots,n_r}{\epsilon_1,\cdots,\epsilon_r}
               =\sum_{0<k_{1}<\cdots<k_{r}} \frac{\epsilon_{1}^{k_{1}}\cdots\epsilon_{r}^{k_{r}}}{k_{1}^{n_{1}}\cdots k_{r}^{n_{r}}}, \epsilon_{i}\in \{1,-1\},(n_r,\epsilon_r)\neq(1,1)
,$$ and we call $N=n_1+\cdots+n_r$ and $r$ its weight and depth respectively.
 Denote by $\mathcal{D}_r\mathcal{Z}^{2}$ the subspace of $\mathcal{Z}^{2}$ spanned by elements
of depth $\leq r$,
then we have $\zeta^{o} (n_{1},n_{2})\in \mathcal{D}_2\mathcal{Z}^{2}$ according to the above discussion.

Similarly we have
\[
\begin{split}
\zeta^{o}(n_1,n_2,\cdots,n_r)&=\sum_{0<k_1<k_2<...<k_r,\mathrm{odd}}\frac{1}{k_1^{n_1}k_2^{n_2}\cdots k_r^{n_r}}\\
                                          &=\frac{1}{2^r}\sum_{\epsilon_i\in\{\pm 1\},1\leq i\leq r}\epsilon_1\epsilon_2\cdots \epsilon_r\zeta \binom{n_1,\cdots,n_r}{\epsilon_1,\cdots,\epsilon_r}.
\end{split}
\]
We call $\zeta^{o} (n_{1},n_{2},\cdots,n_r)$  a sum odd multiple zeta value.

In the case of even weight, Kaneko and Tasaka  \cite{ref7} found the following result:

\begin{Thm} (Kaneko, Tasaka)\label{dim}
For any even integer $k\geq4$, denote by $\mathcal{S}_{k}(2)$ the space of cusp forms for $\Gamma_{0}(2)=\{\gamma\in SL_{2}(\mathbb{Z}); \gamma=\left(\begin{array}{cc} a & b\\ c & d\end{array}\right), c \equiv 0\ mod \ 2\}$ of weight $k$, then we have
$$\mathrm{dim}_{\mathbb{Q}} \langle\zeta^{o}(2r,k-2r); 1\leq r\leq k/2-1\rangle_{\mathbb{Q}} \leq k/2-1-\mathrm{dim}_{\mathbb{C}}\ \mathcal{S}_{k}(2).$$
\end{Thm}

Besides, for $k$ even they also conjectured that elements
$$\zeta^{o}(2r-1,k-2r+1), 1\leq r\leq k/2-1$$
are $\mathbb{Q}$-linear independent, and each element  $\zeta^{o}(r,k-r), 1\leq r\leq k-2$  can be written as a $\mathbb{Q}$-linear combination of
$$\zeta^{o}(2r-1,k-2r+1), 1\leq r\leq k/2-1,\zeta(k).$$

 In this paper we will reinterpret  Kaneko and Tasaka results on the motivic level.

There are weight grading and depth filtration structures on $\mathcal{H}$ which are compatible with the usual weight and depth structures on classical multiple zeta values relative to $\mu_2$. Denote by $gr_r^{\mathcal{D}}\mathcal{H}_N$ the weight $N$ depth $r$ part of the depth-graded motivic multiple zeta values, and
\[
gr_{\bullet}\mathcal{H}=\mathbb{Q}\oplus\bigoplus_{r\geq 1}\mathcal{D}_r\mathcal{H}/\mathcal{D}_{r-1}\mathcal{H}.
\]

Denote by $\zeta^{o,\mathfrak{m}} (n_{1},n_{2})$ the motivic sum odd double zeta value, about which we will introduce later. Let $\mathcal{P}_{ev}^{o,\mathfrak{m}}$  be the space generated  by the natural images of
$$\{\zeta^{o,\mathfrak{m}} (n_{1},n_{2}), n_{1}, n_{2}\geq 2,\mathrm{even}\}$$
in $gr_2^{\mathcal{D}}\mathcal{H}$.

Our first main result (in a rough version) is the following:
\begin{Thm}\label{1.3}
There is an exact sequence with respect to sum odd motivic double zeta values
$$0 \rightarrow \mathcal{P}_{ev}^{o,\mathfrak{m}} \rightarrow (gr_{1}^{\mathcal{D}}\mathcal{H}^{odd}\otimes gr_{1}^{\mathcal{D}}\mathcal{H}^{odd})^b\rightarrow \mathbb{P}(\Gamma_{0}(2))^{\vee} \rightarrow 0,$$
where $\mathcal{H}$ is the algebra of motivic multiple zeta values relative to $\mu_2$.
\end{Thm}

Details of the above notations will be introduced in Section  \ref{periodp} and Section \ref{dtwo}. Theorem
\ref{1.3} gives an description of the space of motivic  sum odd multiple zeta values of the form 
\[
\zeta^{\mathfrak{m},o}(n_1,n_2),n_1,n_2\geq 2, \mathrm{even},
\]
and from it we recover Kaneko and Tasaka's Theorem \ref{dim} immediately.

 We can also discuss the case of odd $n_{1}, n_{2}$, and obtain the following theorem:

\begin{Thm}\label{basis}
$(i)$ For an even integer $N\geq 4$,  the set of  the images of elements
$$\{\zeta^{o,\mathfrak{m}} (n_1,n_2), n_1+n_2=N, 1\leq n_i\leq N-1,\mathrm{odd}\}$$ in $gr_2^{\mathcal{D}}\mathcal{H}$ is  a basis for $gr_2^{\mathcal{D}}\mathcal{H}_N$.\\
$(ii)$ For an odd integer $N\geq 5$, the set of the images of  elements
\[
\{\zeta^{o,\mathfrak{m}}(n_1,n_2,n_3),n_1+n_2+n_3=N,1\leq n_i\leq N-2,\mathrm{odd}\}
\]
in $gr_3^{\mathcal{D}}\mathcal{H}$ is a basis for $gr_3^{\mathcal{D}}\mathcal{H}_N$.
\end{Thm}

In the above theorem, $\zeta^{o,\mathfrak{m}}(n_1,n_2,n_3)$ means the motivic version of the sum odd multiple zeta value $\zeta^{o}(n_1,n_2,n_3)$. Its definition will be given in Section \ref{mtwo}.

From the explicit calculations in the proof of Theorem \ref{1.3} and Theorem \ref{basis}, we also obtain  the following theorem, which was conjectured by Kaneko and Tasaka  \cite{ref7} in Section 3.2, Remark 2.
\begin{Thm}\label{eee}
$(i)$ For an even integer $N\geq 4$, the space $$\langle\zeta^{o}(r, N-r); 1\leq r\leq N-2\rangle_{\mathbb{Q}}$$ is spanned by
$$\{\zeta(N),\zeta^{o}(r, N-r); 1\leq r\leq N-3,\mathrm{odd}\}.$$\\
$(ii)$
For an even integer $N\geq6$, we have
$$\langle\zeta(n_{1},n_{2}); n_{1}+n_{2}=N, n_{2}\geq 2\rangle_{\mathbb{Q}}\subseteq \langle\zeta^{o}(n_{1},n_{2}); n_{1}+n_{2}=N, n_{2} \geq 2\rangle_{\mathbb{Q}},$$
	$$\langle\zeta^{o}(n_{1},n_{2}); n_{1}+n_{2}=N, n_{i}\ even\rangle_{\mathbb{Q}}\subseteq \langle\zeta(n_{1},n_{2}); n_{1}+n_{2}=N, n_{2}\geq 2\rangle_{\mathbb{Q}}.$$
\end{Thm}

We can also give some information in higher depth cases, in the case of depth $3$, we have:
\begin{Thm}\label{deth}
 For a given odd integer $N\geq 5$, and $n_1+n_2+n_3=N, n_1,n_2\geq 1,n_3\geq 2$, the element  $$\zeta^{o}(n_1, n_2, n_3)$$
can be written as a $\mathbb{Q}$-linear combination of $$\zeta^{o}(m_1, m_2, m_3), m_1+m_2+m_3=N, m_1,m_2\geq 1,m_3\geq 3,m_i\,\mathrm{odd}$$ and lower depth multiple zeta values relative to $\mu_2$.
\end{Thm}

It seems that Theorem \ref{basis} should also be true for higher depth. We calculate the depth-graded motivic Galois action for sum odd motivic multiple zeta values explicitly in higher depth. We show that if a special matrix is invertible, then we can prove the higher depth analogue of Theorem \ref{basis}.

Our paper is in the following several parts. In Subsection \ref{mtm}, we introduce mixed Tate motives over $\mathbb{Z}[\frac{1}{2}]$. In Subsection \ref{mmzv}, we  introduce motivic
multiple zeta values relative to $\mu_2$, which was defined by Glanois \cite{ref5}.
We consider the motivic Galois action and show the way to do the calculation
in Subsection \ref{partial}. Then we give a brief
introduction of period polynomials in Section \ref{periodp}. The proofs of our main results will
be given in Section \ref{dtwo} and Section \ref{hdc}.

\section{Motivic multiple zeta values relative to $\mu_2$}\label{mtwo}

As we said in the introduction, the motivic multiple zeta values relative to $\mu_N$ are the generalization of Brown's motivic multiple zeta values. In this section, we only define them in the case of $N=2$. The main references of this section are  \cite{ref4}, \cite{ref5} and \cite{ref6}.

\subsection{Mixed Tate motives over $\mathbb{Z}[\frac{1}{2}]$}\label{mtm}
Consider the category of mixed Tate motives over $\mathbb{Z}[\frac{1}2]$, we denote it by $\mathcal{MT}_{2}$. It is a Tannakian category with the natural fiber functor
$$\omega: \mathcal{MT}_{2}\rightarrow Vec_{\mathbb{Q}}; M\mapsto \oplus \omega_{r}(M),$$
where
$$\omega_{r}(M)=Hom_{\mathcal{MT}_2}(\mathbb{Q}(r), gr_{-2r}^{\omega}(M)).$$

Let $\mathcal{G}^{\mathcal{MT}_{2}}$ be the Tannakian fundamental group (the motivic Galois group) of $\mathcal{MT}_{2}$ with respect to this fiber functor $\omega$, and $\mathcal{U}^{\mathcal{MT}_{2}}$ be the pro-unipotent radical of $\mathcal{G}^{\mathcal{MT}_{2}}$. We have
$$\mathcal{G}^{\mathcal{MT}_{2}}\cong\mathbb{G}_{m}\ltimes \mathcal{U}^{\mathcal{MT}_{2}}.$$

The extension group $\mathrm{Ext}_{\mathcal{MT}_{2}}^{1}(\mathbb{Q}(0), \mathbb{Q}(n))$ is non-trivial only when $n\geq 1, \mathrm{odd}$ and:
$$\mathrm{Ext}_{\mathcal{MT}_{2}}^{1}(\mathbb{Q}(0), \mathbb{Q}(n))\cong\mathbb{Q},\ n\geq 1,\mathrm{ odd},$$
$$\mathrm{Ext}_{\mathcal{MT}_{2}}^{2}(\mathbb{Q}(0), \mathbb{Q}(n))=0,\ \forall n.$$

By the standard argument, there is a set of symbols $\{f_{2n+1}; n\geq 0\}$ such that
$$\mathcal{O}(\mathcal{U}^{\mathcal{MT}_{2}})\cong\mathbb{Q}\langle f_{1},f_{3},\cdots,f_{2n+1},\cdots\rangle,$$
where $\mathbb{Q}\langle f_{1},f_{3},\cdots,f_{2n+1},\cdots\rangle$ means the non-commutative polynomial ring with variables $f_1,f_3,\cdots,f_{2n+1},\cdots.$

Let $\mathfrak{g}$ be the Lie algebra of $\mathcal{U}^{\mathcal{MT}_{2}}$, then $\mathfrak{g}= (m/m^{2})^{\vee}$, where $m\subseteq \mathcal{O}(\mathcal{U}^{\mathcal{MT}_{2}})$ is the maximal ideal. It is a free Lie algebra with a set of generators $\{\sigma_{2n+1}; n\geq 0\}$.

Denote by ${}_{0}\Pi_{1}=\pi_{1}^{dR}(\mathbb{P}^{1}-\{0,1,-1,\infty\}, \overrightarrow{1_{0}},\overrightarrow{-1_{1}})$ the  motivic torsor of paths from $0$ to $1$ on $\mathbb{P}^{1}-\{0,\pm1,\infty\}$, with tangential base point given by the tangent vectors $1$ at $0$ and $-1$ at $1$. It is a functor. For any $\mathbb{Q}$-algebra $R$, denote by $R\langle\langle e_0,e_{-1},e_{1}\rangle\rangle$ the non-commutative $R$-coefficients formal power series in $e_0,e_{-1}, e_1$ and  $$\bigtriangleup:R\langle\langle e_0,e_{-1},e_1\rangle\rangle\to R\langle\langle e_0,e_{-1},e_1\rangle\rangle\otimes_{R} R\langle\langle e_0,e_{-1},e_1\rangle\rangle$$ the co-product on $R\langle\langle e_{0},e_{-1},e_{1}\rangle\rangle$  satisfying $\bigtriangleup e_i=e_i\otimes 1+1\otimes e_i$ for $i\in \{0,\pm 1\}$. Let $R\langle\langle e_0,e_{-1},e_1\rangle \rangle^{\times}$ be the set of non-zero elements of $R\langle\langle e_0,e_{-1},e_1\rangle\rangle$,  we have
$${}_0\Pi_{1}(R)=\{S\in R\langle\langle e_{0},e_{-1},e_{1}\rangle\rangle^{\times};\ \bigtriangleup S=S\otimes S\},$$
i.e. ${}_0\Pi_{1}(R)$ is the set of group-like elements in $R\langle\langle e_{0},e_{-1},e_{1}\rangle\rangle^{\times}$.

Denote by $e^{i}$ be the natural dual of $e_{i}$ for $i\in \{0,1,-1\}$, the affine ring of regular functions of ${}_0\Pi_{1}$ is the graded algebra with the shuffle product
$$\mathcal{O}({}_0\Pi_{1})\cong \mathbb{Q}\langle e^{0},e^{1},e^{-1}\rangle.
$$
The symbol ${}_0 1_1$ is the point ${}_0 1_1:\mathrm{Spec}\;\mathbb{Q}\to {}_0\Pi_1$ whose function ring homomorphism maps every non empty words in $e^0,e^1,e^{-1}$ to $0$.

More generally, for $x,y\in\{0,\pm 1\}$, denote by ${}_x\Pi_y$ the motivic fundamental groupoid from the tangential point at $x$ to the tangential point at $y$.

Let $V$ be the automorphism subgroup of the motivic fundamental groupoid (all basepoints are tangential points at $\{0,\pm 1\}$) of $\mathbb{P}^1-\{0,\pm 1,\infty\}$ satisfying the following properties:\\
(i) Elements of $V$ are compatible with the composition law on the motivic groupoid of $\mathbb{P}^1-\{0,\pm 1,\infty\}$;\\
(ii) Elements of $V$ fix $\mathrm{exp}(e_i)\in {}_i\Pi_i$ for $i\in \{0,\pm1\}$;\\
(iii) Elements of $V$ are equivariant with the $\{\pm 1\}$-action on the motivic groupoid.

By Proposition $5.11$ in \cite{ref4}, the map
\[
\xi:V\to {}_0\Pi_1,
a\mapsto a({}_0 1_1)
\]
is an isomorphism of schemes and
\[
\mathrm{Lie}\,V=(\mathbb{L}(e_0,e_1,e_{-1}),\{\,,\,\}),
\]
where $\mathbb{L}(e_{0},e_{1},e_{-1})$ is the free Lie algebra generated by the three symbols $e_{0},e_{1},e_{-1}$,  and $\{\,,\,\}$ denotes the Ihara Lie bracket on $\mathbb{L}(e_0,e_1,e_{-1})$.

The action of $\mathcal{U}^{\mathcal{MT}_2}$ on ${}_x\Pi_y, x,y\in\{0,\pm 1\}$ factors through $V$.  So there is a natural Lie algebra homomorphism:
$$i: \mathfrak{g}\to \mathrm{Lie}\, V= (\mathbb{L}(e_{0},e_{1},e_{-1}),\{\ ,\ \}).$$

By the main results of Deligne \cite{ref4}, the map $i$ is injective.

 For any element $w$ in $\mathbb{L}(e_{0},e_{1},e_{-1})$, let $\mathrm{depth}(w)$ be the smallest number of  total occurrences of $e_{1}$ and $e_{-1}$ in $w$. It induces a depth filtration $\mathcal{D}$ on $\mathbb{L}(e_{0},e_{1},e_{-1})$ as follows:
$$\mathcal{D}^r\mathbb{L}(e_0,e_1,e_{-1})=\{w \in \mathbb{L}(e_0,e_1,e_{-1}); \mathrm{depth}(w)\geq r\}.$$

According to \cite{ref4},  the map $i$ satisfies:
$$i(\sigma_{1})=e_{-1},$$
$$i(\sigma_{2n+1})=(1-2^{2n})\mathrm{ad}(e_{0})^{2n}e_{-1}+2^{2n}\mathrm{ad}(e_{0})^{2n}e_{1}+\mathrm{HDT},$$
where $\mathrm{HDT}$ means the higher depth term.

The motivic Lie algebra $\mathfrak{g}$ has an induced depth filtration $\mathcal{D}^r\mathfrak{g}$ from the injective map $i$. Since Ihara bracket is compatible with the depth filtration, we know that the depth-graded space
\[
\mathfrak{dg}=\oplus_{r\geq 1}\mathcal{D}^r\mathfrak{g}/\mathcal{D}^{r+1}\mathfrak{g}
\]
 is a Lie algebra with the induced Ihara Bracket. Furthermore, from the main results of  Deligne \cite{ref4}, $\mathfrak{dg}$ is a free Lie algebra with generators
$$\overline {i(\sigma_{1})}=e_{-1},$$
$$\overline{i(\sigma_{2n+1})}=(1-2^{2n})\mathrm{ad}(e_{0})^{2n}e_{-1}+2^{2n}\mathrm{ad}(e_{0})^{2n}e_{1}$$ in the depth one part.

We will use them in the style of Lie polynomial in $\mathbb{Q}\langle e_0,e_1,e_{-1}\rangle$ rather than Lie words in the rest of this paper for convenience:
$$i(\sigma_{2n+1})=(1-2^{2n})\sum_{r=0}^{2n}(-1)^r{2n \choose r}e_{0}^{2n-r}e_{-1}e_{0}^{r}+2^{2n}\sum_{r=0}^{2n}(-1)^r{2n \choose r}e_{0}^{2n-r}e_{1}e_{0}^{r}+\mathrm{HDT}.$$

\subsection{Motivic multiple zeta values}\label{mmzv}
Similar to Brown's work, Glanois \cite{ref5} defined motivic iterated integral $I^{\mathfrak{m}}$ and motivic multiple zeta values
$\zeta^{\mathfrak{m}}\binom{x_1,x_2,...,x_p}{\epsilon_1,\epsilon_2,...,\epsilon_p}$, $\epsilon_{i}\in \mu_N$ relative to the set of $N^{th}$ roots of unity $\mu_N$.
We denote by $\mathcal{H}=\mathcal{H}^{2}$ the $\mathbb{Q}$-vector space of motivic multiple zeta values relative to $\mu_2=\{1,-1\}$. Here we only give the definition in the case of $N=2$.

Now we construct the  map
\[
{dch}:\mathbb{Q}\langle e^0,e^1,e^{-1}\rangle\to \mathbb{R}
\]
for words $u_{i}\in \{e^{0},e^{1},e^{-1}\}, i=1,\cdots,k$ satisfying $u_{1}\neq e^{0}, u_{k}\neq e^{1}$, define
$${dch}(u_{1}\cdots u_{k})=\int_{0<t_{1}<\cdots <t_{k}<1}\omega_{u_{1}}(t_{1})\cdots \omega_{u_{k}}(t_{k}),$$
where $\omega_{e^{0}}(t)=dt/t, \omega_{e^{i}}(t)=dt/(i-t), i\in \{1,-1\}$.
In general, we know that
$$\int_{\varepsilon<t_{1}<\cdots <t_{k}<1-\eta}\omega_{u_{1}}(t_{1})\cdots \omega_{u_{k}}(t_{k})=P(log(\varepsilon), log(\eta))+O(\mathrm{sup}(\varepsilon|log(\epsilon)|^{A}+\eta|log(\eta)|^{B})),$$
where $P$ is a polynomial. For a general word sequence $u_1\cdots u_k$,  define $${dch}(u_1,...,u_k)=P(0,0).$$

By the shuffle product of iterated integral, ${dch}$ is a $\mathbb{Q}$-algebra homomorphism.
Denote by $\mathcal{I}$ the largest graded sub-ideal of $\mathrm{Ker}\,{dch}$ which is stable under the action of $\mathcal{G}^{\mathcal{MT}_2}$.
The motivic multiple zeta algebra $\mathcal{H}$ is  $\mathcal{O}({}_0\Pi_{1})/\mathcal{I}$.

Denote by $\mathcal{I}^{\mathfrak{m}}$ the natural quotient map
\[
\mathcal{I}^{\mathfrak{m}}:\mathcal{O}({}_0\Pi_1)=\mathbb{Q}\langle e^0, e^1, e^{-1}\rangle\to \mathcal{H}
\]
and $per$ the map $per:\mathcal{H}\to \mathbb{R}$  satisfying $per\circ \mathcal{I}^{\mathfrak{m}}={dch}$.

The motivic multiple zeta value $\zeta^{\mathfrak{m}}\binom{x_1,x_2...,x_p}{\epsilon_1,\epsilon_2,...,\epsilon_p}$ is
\[
\mathcal{I}^{\mathfrak{m}}(e^{(\epsilon_1\cdots \epsilon_p)^{-1}}(e^0)^{x_1-1}e^{(\epsilon_2\cdots \epsilon_p)^{-1}} (e^0)^{x_2-1}\cdots e^{(\epsilon_p)^{-1}}  (e^0)^{x_p-1}             ).
\]
It's obvious to check that  \[
per\left(\zeta^{\mathfrak{m}}\binom{x_1,...,x_p}{\epsilon_1,...,\epsilon_p}\right)=\zeta\binom{x_1,...,x_p}{\epsilon_1,...,\epsilon_p}.
\]
Define $\zeta^{o,\mathfrak{m}}(n_{1},\cdots,n_{r})$ as
\[
\begin{split}
&\zeta^{o,\mathfrak{m}}(n_{1},\cdots,n_{r})
=\frac{1}{2^{r}}\sum_{\epsilon_i\in \{\pm1\},1\leq i\leq r}\epsilon_1\cdots \epsilon_r\zeta^{\mathfrak{m}}\binom{n_1,\cdots,n_r}{\epsilon_1,\cdots, \epsilon_r}.
\end{split}
\]
It's clear that the image of $\zeta^{o,\mathfrak{m}}(n_1,...,n_r)$ under the period map $per$ is the sum odd multiple zeta values $\zeta^{o}(n_1,...,n_r)$.

In $\mathcal{O}({}_0\Pi_1)=\mathbb{Q}\langle e^0,e^1,e^{-1}\rangle$, for any word $u_1\cdots u_k, u_i\in \{e^0,e^1,e^{-1}\}$, $k$ is called its weight and the total number of occurrences of $e^1$ and $e^{-1}$ is called its depth. Denote by $\mathcal{D}_r\mathbb{Q}\langle e^0, e^1, e^{-1}\rangle$ the subspace spanned by elements of depth $\leq r$.

 Since the depth filtration on $\mathcal{O}({}_0\Pi_1)$ is motivic \cite{ref3}, it induces a natural depth filtration on $\mathcal{H}$. The depth filtration on $\mathcal{H}$ is compatible with the depth filtration on $\mathcal{Z}^2$ through the map $per$.

Denote by $gr_r^{\mathcal{D}}\mathcal{H}=\mathcal{D}_r\mathcal{H}/\mathcal{D}_{r-1}\mathcal{H}$, Deligne and Goncharov proved the following formula in the case of depth $1$, as a consequence of this formula we can find a basis of $gr_{1}^{\mathcal{D}}\mathcal{H}$:
\begin{lem} \label{2.1}(Deligne-Goncharov)
We have the distribution formula
$$\zeta^{\mathfrak{m}}\binom{n}{-1}
= (2^{-n+1}-1)\zeta^{m}\binom{n}{1},\ \forall \,n\geq 2.
$$
\end{lem}

\begin{lem}\label{2.2}(Deligne-Goncharov)
There is a basis of $gr_{1}^{\mathcal{D}}\mathcal{H}$: $\{\zeta^{\mathfrak{m}}\binom{r}{-1},\ r\geq 1\ odd
\}$.
\end{lem}

\begin{rem}
We will always write $\zeta^{\mathfrak{m}}\binom{n_1,n_2}{-1,1}
$ to be $\zeta^{\mathfrak{m}}(\overline{n_{1}},n_{2})$, similarly $\zeta^{\mathfrak{m}}(n_{1},\overline{n_{2}})$, $\zeta^{\mathfrak{m}}(\overline{n_{1}},\overline{n_{2}})$, $\zeta^{\mathfrak{m}}(\overline{k})$ for convenience.
\end{rem}

\subsection{Motivic Galois action}\label{partial}

In this subsection we explain how to calculate the depth-graded version motivic Galois action of the motivic Lie algebra of $\mathcal{MT}(\mathbb{Z}[\frac{1}{2}])$ on the motivic multiple zeta values relative to $\mu_2$. Then we give the definition of the map $\partial$ and deduce its injectivity from the results of Brown \cite{ref2}, Deligne \cite{ref4} and Glanois \cite{ref5}.

Since the expression of $i(\sigma_{2n+1})$ in $(\mathbb{L}(e_0,e_1,e_{-1}),\{\,,\,\})$ has canonical depth one part, $\sigma_{2n+1}$ in $\mathfrak{g}= \mathrm{Lie}\,\,\mathcal{U}^{\mathcal{MT}_{2}}$ induces a well defined derivation $$\partial_{2n+1}:gr_{r}^{\mathfrak{D}}\mathcal{H}\rightarrow gr_{r-1}^{\mathfrak{D}}\mathcal{H}.$$
 In this section we will show how to calculate the map $\partial_{2n+1}$ explicitly.

Since $\mathcal{O}({}_0\Pi_{1})$ is an ind-object in the category $\mathcal{MT}_{2}$, there is an action of the motivic Lie algebra
$$\mathfrak{g}\times \mathcal{O}({}_0\Pi_{1})\rightarrow \mathcal{O}({}_0\Pi_{1}).$$

Denote by $\mathfrak{h}=\mathrm{Lie}\, V=(\mathbb{L}(e_0,e_1,e_{-1}),\{\,,\,\})$. The action of $\mathfrak{g}$ on $\mathcal{O}({}_0\Pi_1)$ factors through the action of $\mathfrak{h}$ on $\mathcal{O}({}_0\Pi_1)$.

Denote by $\mathcal{U}\mathfrak{h}$ the universal enveloping algebra of $\mathfrak{h}$.  Then
\[
\mathcal{U}\mathfrak{h}=(\mathbb{Q}\langle e_0,e_1,e_{-1}\rangle,\circ),\]
 where $\circ$ denotes the new product on $\mathbb{Q}\langle e_0,e_1,e_{-1}\rangle$ transformed from the natural concatenation product on $\mathcal{U}\mathfrak{h}$.

 The product $\circ$  is difficult to calculate in general, but by the same reason as Proposition $2.2$ in \cite{ref1}, for any $a\in \mathfrak{h}$, any non-empty word $w$ in $e_0, e_1,e_{-1}$, and any $n\geq 0$, we have
$$a\circ (e_{0}^{n}e_{i}w)=e_{0}^{n}(a\circ e_{i})w+e_{0}^{n}e_{i}(a\circ w), i\in \{1,-1\},$$
where
$$a\circ e_{0}^n=e_0^n a,\ a\circ e_{1}=ae_{1}+e_{1}a^{*},\ a\circ e_{-1}=([-1]a)e_{-1}+e_{-1}([-1]a)^{*},$$
$$[-1](e_{0}^{a}e_{1}^{b}e_{-1}^{c}\cdots)=e_{0}^{a}e_{-1}^{b}e_{1}^{c}\cdots,\ (a_{1}\cdots a_{n})^{*}=(-1)^{n}(a_{n}\cdots a_{1}), a_i\in \{e_0, e_1, e_{-1} \}.$$

From the correspondence between unipotent algebraic group and nilpotent Lie algebra, we know that for any $a\in\mathfrak{h}$, the natural action of $a$ on $\mathcal{O}({}_0\Pi_1)$:
\[
\mathcal{O}({}_0\Pi_1)=\mathbb{Q}\langle e^0,e^1,e^{-1}\rangle\xrightarrow{a} \mathcal{O}({}_0\Pi_1)=\mathbb{Q}\langle e^0,e^1,e^{-1}\rangle,\]
\[ x\mapsto a(x)
\]
is dual to the following action of $a$ on $\mathcal{U}\mathfrak{h}$:
\[
\mathcal{U}\mathfrak{h}=\mathbb{Q}\langle e_0,e_1,e_{-1}\rangle\xrightarrow{a} \mathcal{U}\mathfrak{h}=\mathbb{Q}\langle e_0,e_1,e_{-1}\rangle,
\]
\[
y\mapsto a\circ y.
\]

By the definition of $\mathcal{H}$ and $\partial_{2n+1}$, we have the following commutative diagram:
 \[
 \xymatrix{
   gr_r^{\mathcal{D}}\mathbb{Q}\langle e^0,e^1,e^{-1}\rangle \ar@{->>}[d] \ar[r]^{\overline{\partial_{2n+1}}} & gr_{r-1}^{\mathcal{D}}\mathbb{Q}\langle e^0,e^1,e^{-1}\rangle \ar@{->>}[d] \\
 gr_r^{\mathcal{D}}\mathcal{H} \ar[r]^{\partial_{2n+1}} & gr_{r-1}^{\mathcal{D}}\mathcal{H} , }
 \]
 where $\overline{\partial_{2n+1}}$ is  the depth-graded version of the action of $i(\sigma_{2n+1})$ on $\mathbb{Q}\langle e^0,e^1,e^{-1}\rangle$.
Thus in order to write out the maps $\overline{\partial_{2n+1}}$ and $\partial_{2n+1}$ clearly, we need to compute the action $\circ: \mathfrak{h}\times \mathcal{U}\mathfrak{h}\rightarrow \mathcal{U}\mathfrak{h}$ first.

There is a well-defined map
$$\partial:gr_r^{\mathcal{D}}\mathcal{H} \rightarrow gr_{1}^{\mathcal{D}}\mathcal{H}^{odd}\otimes gr_{r-1}^{\mathcal{D}}\mathcal{H};\ \partial=\sum_{n\geq0}\zeta^{\mathfrak{m}}(\overline{2n+1})\otimes \partial_{2n+1}.$$
The following Proposition is crucial to our analysis.
\begin{prop}\label{inj}
For $r\geq 2$, the map $\partial$ is injective.
\end{prop}
\noindent{\bf Proof:}
From the main results of \cite{ref2},\cite{ref5}, it follows that
 $\mathcal{H}\cong \mathcal{O}(\mathcal{U}^{\mathcal{MT}_2}) [t]$ ($t$ is a weight $2$, depth $1$ element with trivial action of $\mathfrak{g}$) as a $\mathfrak{g}$-module. What's more $t^n,n\geq 1$ are all depth $1$ elements.

 So we have
 \[
 gr_r\mathcal{H}\cong gr_r\mathcal{O}(\mathcal{U}^{\mathcal{MT}_2})\oplus \bigoplus_{n\geq 1}gr_{r-1}\mathcal{O}(\mathcal{U}^{\mathcal{MT}_2})t^n.
 \]
 It suffices to prove that $\partial|_{gr_r\mathcal{O}(\mathcal{U}^{\mathcal{MT}_2})}$ is injective.
 Since the depth-graded motivic Lie algebra $\mathfrak{dg}$ is a free Lie algebra with generators in the depth one part by the main results of \cite{ref4}. By the correspondence between nilpotent Lie algebra and unipotent algebraic group, $\partial|_{gr_r\mathcal{O}(\mathcal{U}^{\mathcal{MT}_2})}$ is injective.
$\hfill\Box$\\

\begin{rem}\label{com}
Proposition \ref{inj} is not true for Brown's original motivic multiple zeta values. Since in that case, the depth-graded motivic Lie algebra of $\mathcal{MT}(\mathbb{Z})$ is not a free Lie algebra and it has generators in higher depth part. See \cite{ref1}, \cite{en}, \cite{li} for some conjectural descriptions of the depth-graded motivic Lie algebra of $\mathcal{MT}(\mathbb{Z})$.
\end{rem}

\section{Period Polynomials}\label{periodp}

In this section, we review the theory of  period polynomials, and define $\mathbb{P}(\Gamma_{0}(2))^{\vee}$ in  Theorem \ref{1.3}. The main reference is Kaneko, Tasaka \cite{ref7}.

As we know, $\Gamma_{0}(2)=\{\gamma\in SL_{2}(\mathbb{Z}); \gamma \equiv 0\ mod \ 2\}$ is generated by two elements

$$T=\left(
      \begin{array}{cc}
        1 & 1 \\
        0 & 1 \\
      \end{array}
    \right)
, M=\left(
   \begin{array}{cc}
     -1 & -1 \\
     2 & 1 \\
   \end{array}
 \right).
$$

For a positive even integer $k$, denote by $V_{k}$ the space of homogeneous polynomials with two variables $X,\ Y$ of degree $k-2$:
$$V_{k}=\{P(X,Y)\in \mathbb{Q}[X, Y]; P(X,Y)=\sum_{i=0}^{k-2} a_{i}X^{i}Y^{k-2-i}\}.$$

The group $\Gamma= \Gamma_{0}(2)$ acts on $V_{k}$ naturally: for any polynomial $P(X,Y)\in V_{k}$ and $\gamma= \left(
                                                           \begin{array}{cc}
                                                             a & b \\
                                                             c & d \\
                                                           \end{array}
                                                         \right) \in \Gamma_{0}(2)$, $$\gamma \circ P(X,Y)= P(aX+bY, cX+dY),$$ we write this action as $P(X, Y)|\gamma$ for convenience. Consider the subspace $W_{k}$ of $V_{k}$ as follows:
$$W_{k}=\{P(X,Y)\in V_{k}[X, Y]; P|(1-T)(1+M)=0\}.$$

Denote by $\mathcal{S}_k(2)$ the space of cusp forms of weight $k$ for $\Gamma_{0}(2)$. For $f\in\mathcal{S}_k(2)$, the period polynomial $r_{f}(X, Y)$ of $f$ is given by
$$r_{f}(X,Y)=\int_{0}^{\infty}f(z)(X-zY)^{k-2}dz.$$
It can be shown that
$$r_{f}(X,Y) \in W_{k}\otimes \mathbb{C}.$$

Now we decompose $W_{k}$ into two parts. Put $\varepsilon=\left(
                                                            \begin{array}{cc}
                                                              -1 & 0 \\
                                                              0 & 1 \\
                                                            \end{array}
                                                          \right)
$, it is obvious to see that $P|(1\pm \varepsilon)\in W_{k}$ for any $P\in W_{k}$, thus we have the direct sum decomposition:
$$W_{k}=W_{k}^{+}\oplus W_{k}^{-},$$
where $W_{k}^{+}$ (resp. $W_{k}^{-}$) is the even (resp. odd) part of $W_{k}$:
$$W_{k}^{\pm}= \{P\in V_{k}; P|\varepsilon= \pm P, P|(1-T)(1+M)=0\}.$$

For $f\in S_{k}(2)$, denote by $r_{f}^{\pm}$ the even and odd parts of the map $r_{f}$:
$$r^{\pm}: S_{k}(2)\rightarrow W_{k}^{\pm}\otimes \mathbb{C}; f\mapsto r_{f}^{\pm}(X,Y).$$

We can decompose $W_{k}^{+}$ further as
$$W_{k}^{+}=\mathbb{Q}\cdot X^{k-2}\oplus W_{k}^{+,0}\oplus \mathbb{Q}\cdot Y^{k-2},$$
where
$$W_{k}^{+,0}=\{P(X,Y)=\sum_{i=2\ even}^{k-4}a_{i}X^{i}Y^{k-2-i}\in V_{k};P|(1-T)(1+M)=0\}.$$

Kaneko and Tasaka \cite{ref7} proved the following two propositions which describe the structure of $W_{k}^{+,0}$:

\begin{prop}
For any even integer $k$, there are two isomorphisms of vector spaces

$$r^{+}:S_{k}(2)\rightarrow W_{k}^{+,0}\otimes \mathbb{C}\ and\ r^{-}:S_{k}(2)\rightarrow W_{k}^{-}\otimes \mathbb{C}.$$
\end{prop}

\begin{prop}\label{3.2}
For any even integer $n,k$, denote by ${n \choose k}$ the combination number, the space $W_{k}^{+,0}$ is in the following form:
$$\{\sum_{i=2\ evev}^{k-4}a_{i}X^{i}Y^{k-2-i};\sum_{i=2\ even}^{k-4}\left({i \choose j}-{i \choose k-2-j}\right)a_{k-2-j}=0, 1\leq j \leq k-3, odd\}.$$
\end{prop}

\begin{rem}
In this paper,  let ${n \choose k}=0$ when $k>n$ or $k<0$.
\end{rem}

Denote by  $\mathbb{P}(\Gamma_{0}(2))=\bigoplus_{k} W_{k}^{+,0}$. Then $\mathbb{P}(\Gamma_{0}(2))^{\vee}$ in Theorem \ref{1.3} is the compact dual of $\mathbb{P}(\Gamma_{0}(2))$.

\section{The depth two case}\label{dtwo}
In this section we calculate the map $$\partial:gr_r^{\mathcal{D}}\mathcal{H} \rightarrow gr_{1}^{\mathcal{D}}\mathcal{H}^{odd}\otimes gr_{r-1}^{\mathcal{D}}\mathcal{H}$$ in the case  of $r=2$ explicitly. Then we establish a short exact sequence about sum odd motivic double zeta values and we find a basis for the depth-graded motivic double zeta values relative to $\mu_2$ by the explicit expression of the
 map $\partial$ in the case of $r=2$. As an application of our results, we prove Kaneko and Tasaka's conjectures in Remark 2, \cite{ref7}.

\subsection{The calculation in depth two}

The following formulas come from direct calculation.
We write $\overline{i(\sigma_{2n+1})}$ as $\overline{\sigma_{2n+1}}$ for short.
When $n=0,\ \overline{\sigma_{2n+1}}=\overline{\sigma_{1}}$,  $a_1\geq 0$,  we have
$$\overline{\sigma_{1}}\circ(e_{0}^{a_{0}}e_{1}e_{0}^{a_{1}})=
e_{0}^{a_{0}}e_{-1}e_{1}e_{0}^{a_{1}}-e_{0}^{a_{0}}e_{1}e_{-1}e_{0}^{a_{1}}+e_{0}^{a_{0}}e_{1}e_{0}^{a_{1}}e_{-1},$$
$$\overline{\sigma_{1}}\circ(e_{0}^{a_{0}}e_{-1}e_{0}^{a_{1}})=
e_{0}^{a_{0}}e_{1}e_{-1}e_{0}^{a_{1}}-e_{0}^{a_{0}}e_{-1}e_{1}e_{0}^{a_{1}}+e_{0}^{a_{0}}e_{-1}e_{0}^{a_{1}}e_{-1}.$$

When $n>0$,  $a_1\geq 0$, we have
$$\overline{\sigma_{2n+1}}\circ(e_{0}^{a_{0}}e_{1}e_{0}^{a_{1}})=
e_{0}^{a_{0}}(\overline{\sigma_{2n+1}}\circ e_{1})e_{0}^{a_{1}}+e_{0}^{a_{0}}e_{1}e_{0}^{a_{1}}\overline{\sigma_{2n+1}}$$
$$=(1-2^{2n})\sum_{r=0}^{2n}(-1)^{r}{2n \choose r}e_{0}^{a_{0}}e_{0}^{2n-r}e_{-1}e_{0}^{r}e_{1}e_{0}^{a_{1}}$$
$$+2^{2n}\sum_{r=0}^{2n}(-1)^{r}{2n \choose r}e_{0}^{a_{0}}e_{0}^{2n-r}e_{1}e_{0}^{r}e_{1}e_{0}^{a_{1}}$$
$$-(1-2^{2n})\sum_{r=0}^{2n}(-1)^{r}{2n \choose r}e_{0}^{a_{0}}e_{1}e_{0}^{r}e_{-1}e_{0}^{2n-r}e_{0}^{a_{1}}$$
$$-2^{2n}\sum_{r=0}^{2n}(-1)^{r}{2n \choose r}e_{0}^{a_{0}}e_{1}e_{0}^{r}e_{1}e_{0}^{2n-r}e_{0}^{a_{1}}$$
$$+(1-2^{2n})\sum_{r=0}^{2n}(-1)^{r}{2n \choose r}e_{0}^{a_{0}}e_{1}e_{0}^{a_{1}}e_{0}^{2n-r}e_{-1}e_{0}^{r}$$
$$+2^{2n}\sum_{r=0}^{2n}(-1)^{r}{2n \choose r}e_{0}^{a_{0}}e_{1}e_{0}^{a_{1}}e_{0}^{2n-r}e_{1}e_{0}^{r},$$
and
$$\overline{\sigma_{2n+1}}\circ(e_{0}^{a_{0}}e_{-1}e_{0}^{a_{1}})=
e_{0}^{a_{0}}(\overline{\sigma_{2n+1}}\circ e_{-1})e_{0}^{a_{1}}+e_{0}^{a_{0}}e_{-1}e_{0}^{a_{1}}\overline{\sigma_{2n+1}}$$
$$=(1-2^{2n})\sum_{r=0}^{2n}(-1)^{r}{2n \choose r}e_{0}^{a_{0}}e_{0}^{2n-r}e_{1}e_{0}^{r}e_{-1}e_{0}^{a_{1}}$$
$$+2^{2n}\sum_{r=0}^{2n}(-1)^{r}{2n \choose r}e_{0}^{a_{0}}e_{0}^{2n-r}e_{-1}e_{0}^{r}e_{-1}e_{0}^{a_{1}}$$
$$-(1-2^{2n})\sum_{r=0}^{2n}(-1)^{r}{2n \choose r}e_{0}^{a_{0}}e_{-1}e_{0}^{r}e_{1}e_{0}^{2n-r}e_{0}^{a_{1}}$$
$$-2^{2n}\sum_{r=0}^{2n}(-1)^{r}{2n \choose r}e_{0}^{a_{0}}e_{-1}e_{0}^{r}e_{-1}e_{0}^{2n-r}e_{0}^{a_{1}}$$
$$+(1-2^{2n})\sum_{r=0}^{2n}(-1)^{r}{2n \choose r}e_{0}^{a_{0}}e_{-1}e_{0}^{a_{1}}e_{0}^{2n-r}e_{-1}e_{0}^{r}$$
$$+2^{2n}\sum_{r=0}^{2n}(-1)^{r}{2n \choose r}e_{0}^{a_{0}}e_{-1}e_{0}^{a_{1}}e_{0}^{2n-r}e_{1}e_{0}^{r}.$$

By taking dual of the first formula, we have the following result:
\begin{lem}\label{4.1}
For positive even integers $n_{1}, n_{2}$ and $\epsilon_{1}, \epsilon_{2}\in\{1, -1\}$
$$\overline{\partial_{1}}(e^{\epsilon_{1}}(e^{0})^{n_{1}-1}e^{\epsilon_{2}}(e^{0})^{n_{2}-1})=0.$$
\end{lem}
\noindent{\bf Proof:} We calculate the map by taking dual of the action $\overline{\sigma_{2n+1}}$, thus we only need to find the terms $e_{\epsilon_{1}}(e_{0})^{n_{1}-1}e_{\epsilon_{2}}(e_{0})^{n_{2}-1}$ on the right hand side of the equation in the first formula. However, there are no such terms because $n_{1}, n_{2}$ are both even and thus $n_{1}-1,\ n_{2}-1$ are odd, it means that there is at least one $e_{0}$ between $e_{\epsilon_{1}}$ and $e_{\epsilon_{2}}$, and one $e_{0}$ after $e_{\epsilon_{2}}$. It follows that $\partial_{1}(e^{\epsilon_{1}}(e^{0})^{n_{1}-1}e^{\epsilon_{2}}(e^{0})^{n_{2}-1})=0$ for all $n_{1}, n_{2}$ even.
$\hfill\Box$\\

For the same reason, the following lemma holds:

\begin{lem}\label{4.2}
For positive even integers $n_{1}, n_{2}$ and $n>0$, write the word
$${2n \choose n_{1}-1}e^{1}(e^{0})^{n_{1}+n_{2}-2n-2},\ {2n \choose n_{1}-1}e^{-1}(e^{0})^{n_{1}+n_{2}-2n-2},$$
$${2n \choose n_{2}-1}e^{1}(e^{0})^{n_{1}+n_{2}-2n-2},\ {2n \choose n_{2}-1}e^{-1}(e^{0})^{n_{1}+n_{2}-2n-2},$$
as $\Theta_{1}^{n_{1}}$, $\Theta_{-1}^{n_{1}}$, $\Theta_{1}^{n_{2}}$, $\Theta_{-1}^{n_{2}}$ respectively for convenience, we have
$$\overline{\partial_{2n+1}}(e^{1}(e^{0})^{n_{1}-1}e^{1}(e^{0})^{n_{2}-1})=2^{2n}(\Theta_{1}^{n_{1}}-\Theta_{1}^{n_{2}}),$$
$$\overline{\partial_{2n+1}}(e^{1}(e^{0})^{n_{1}-1}e^{-1}(e^{0})^{n_{2}-1})=(1-2^{2n})(\Theta_{1}^{n_{1}}-\Theta_{1}^{n_{2}}),$$
$$\overline{\partial_{2n+1}}(e^{-1}(e^{0})^{n_{1}-1}e^{1}(e^{0})^{n_{2}-1})=(1-2^{2n})\Theta_{-1}^{n_{1}}-2^{2n}\Theta_{-1}^{n_{2}},$$
$$\overline{\partial_{2n+1}}(e^{-1}(e^{0})^{n_{1}-1}e^{-1}(e^{0})^{n_{2}-1})=2^{2n}\Theta_{-1}^{n_{1}}-(1-2^{2n})\Theta_{-1}^{n_{2}}.$$
\end{lem}

It is also useful for us to determine the case that $n_{1}, n_{2}$ are both odd. We use the same argument here and the result is a little different.

\begin{lem}\label{4.3}
For positive odd integers $n_{1}\geq 1, n_{2}\geq 1$, we have

$$ \partial_{1}(e^{i_{1}}(e^{0})^{n_1-1}e^{i_{2}}(e^0)^{n_2-1})$$
$$=\left\{
\begin{array}{lcl}
-e^{1}(e^0)^{n_2-1}+ e^{-1}(e^0)^{n_2-1},    &      & ({i_{1}, i_{2})=(1,-1), n_1=1,n_2\geq 1},\\
e^{1}(e^0)^{n_2-1}- e^{-1}(e^0)^{n_2-1}  ,   &      & ({i_{1}, i_{2})=(-1,1), n_1=1,n_2\geq 1},\\
e^1(e^0)^{n_1-1},&       & ({i_{1}, i_{2})=(1,-1), n_1\geq 3,n_2=1},\\
e^{-1}(e^0)^{n_1-1},&         & ({i_{1}, i_{2})=(-1,-1), n_1\geq 3},n_2=1,\\
0     ,                                          &      & {otherwise}.\\
\end{array} \right. $$

\end{lem}

Define $\delta\binom{m}{n}=1$ if $m=n$, $\delta\binom{m}{n}=0$ if $m\neq n$.
\begin{lem}\label{4.4}
For positive odd integers $n_{1}, n_{2}$, let $\Theta_{1}^{n_{1}}$, $\Theta_{-1}^{n_{1}}$, $\Theta_{1}^{n_{2}}$, $\Theta_{-1}^{n_{2}}$ be as above and $n\geq 1$, we have
$$\overline{\partial_{2n+1}}(e^{1}(e^{0})^{n_{1}-1}e^{1}(e^{0})^{n_{2}-1})=-2^{2n}\left((\Theta_{1}^{n_{1}}-\Theta_{1}^{n_{2}})-\delta\binom{2n}{n_{1}-1}\Theta_{1}^{n_{1}}\right),$$
$$\overline{\partial_{2n+1}}(e^{1}(e^{0})^{n_{1}-1}e^{-1}(e^{0})^{n_{2}-1})=-(1-2^{2n})\left((\Theta_{1}^{n_{1}}-\Theta_{1}^{n_{2}})-\delta\binom{2n}{n_{1}-1}\Theta_{-1}^{n_{1}}\right),$$
$$\overline{\partial_{2n+1}}(e^{-1}(e^{0})^{n_{1}-1}e^{1}(e^{0})^{n_{2}-1})=-(1-2^{2n})\left(\Theta_{-1}^{n_{1}}-\delta\binom{2n}{n_{1}-1}\Theta_{1}^{n_{1}}\right)+2^{2n}\Theta_{-1}^{n_{2}},$$
$$\overline{\partial_{2n+1}}(e^{-1}(e^{0})^{n_{1}-1}e^{-1}(e^{0})^{n_{2}-1})=-2^{2n}\left(\Theta_{-1}^{n_{1}}-\delta\binom{2n}{n_{1}-1}\Theta_{-1}^{n_{1}}\right)+(1-2^{2n})\Theta_{-1}^{n_{2}}.$$
\end{lem}

With the above lemmas, we can calculate the maps $\overline{\partial}$ and $\partial$ in the case of $r=2$.

\subsection{Proofs of the main results}
Now we are ready to state our main results. We have already defined the map $\partial_{2n+1}$ for $n\geq 0$ and the space $\mathcal{P}_{ev}^{o,\mathfrak{m}}$, which is the subspace of $gr_{2}^{\mathcal{D}}\mathcal{H}$ generated by the set of images of  $\{\zeta^{o,\mathfrak{m}}(n_1,n_2)
,\ n_{1}, n_{2}\geq 2,\mathrm{even}\}$.
Define
$$D:gr_{1}^{\mathcal{D}}\mathcal{H}^{odd}\otimes gr_{1}^{\mathcal{D}}\mathcal{H}^{odd}\rightarrow (\mathbb{P}(\Gamma_{0}(2)))^{{}^{\vee}},$$
$$\zeta^{\mathfrak{m}}(\overline{2n_{1}+1})\otimes\zeta^{\mathfrak{m}}(\overline{2n_{2}+1}) \mapsto \frac{2^{2n_{2}}-1}{2^{2n_{2}+1}-1}\cdot v(2n_1+1,2n_2+1), n_1,n_2\geq 0,
$$
where $v(2n_1+1,2n_2+1)$ is a linear functional on $\mathbb{P}(\Gamma_{0}(2))$ satisfying
\[
v(2n_1+1,2n_2+1)(p)=p_{2n_1,2n_2}
\]
for
\[
p=\sum p_{2m_1,2m_2} X^{2m_1}Y^{2m_2}\in \mathbb{P}(\Gamma_0(2))
.\]
\begin{Thm}\label{ES} Denote by $(gr_{1}^{\mathcal{D}}\mathcal{H}^{odd}\otimes gr_{1}^{\mathcal{D}}\mathcal{H}^{odd})^b$ the subspace of $gr_{1}^{\mathcal{D}}\mathcal{H}^{odd}\otimes gr_{1}^{\mathcal{D}}\mathcal{H}^{odd}$ which is generated by $\zeta^{\mathfrak{m}}(\overline{n}_1)\otimes \zeta^{\mathfrak{m}}(\overline{n}_2), n_1,n_2\geq 3$.
Then $$\partial(\mathcal{P}_{ev}^{o,\mathfrak{m}})\subseteq(gr_{1}^{\mathcal{D}}\mathcal{H}^{odd}\otimes gr_{1}^{\mathcal{D}}\mathcal{H}^{odd})^b$$ and  there is an exact sequence:
$$0 \rightarrow \mathcal{P}_{ev}^{o,\mathfrak{m}} \xrightarrow{\widetilde{\partial}} (gr_{1}^{\mathcal{D}}\mathcal{H}^{odd}\otimes gr_{1}^{\mathcal{D}}\mathcal{H}^{odd})^b\xrightarrow{\widetilde{D}} \mathbb{P}(\Gamma_{0}(2))^{{}^{\vee}} \rightarrow 0,$$
where the second map $\widetilde{\partial}$ is induced from $\partial|_{\mathcal{P}_{ev}^{o,\mathfrak{m}}}$ and the third map $\widetilde{D}$ is induced from $D$ defined as above.
\end{Thm}
\noindent{\bf Proof:}  By Lemma \ref{4.1} and Lemma \ref{4.2} it's obvious to check that
\[
\partial(\mathcal{P}_{ev}^{o,\mathfrak{m}})\subseteq(gr_{1}^{\mathcal{D}}\mathcal{H}^{odd}\otimes gr_{1}^{\mathcal{D}}\mathcal{H}^{odd})^b.
\]
The map $\widetilde{\partial}$ is injective by Proposition \ref{inj}.
The surjectivity of $\widetilde{D}$ is trivial. We only need to show  that $\mathrm{Im}\,\widetilde{\partial}=\mathrm{Ker}\,\widetilde{D}$.

The following  diagram is commutative:
$$\begin{array}[c]{ccc}
gr_{2}^{\mathcal{D}}\mathbb{Q}\langle e^{0},e^{1},e^{-1}\rangle&\stackrel{\overline{\partial}}{\rightarrow}&(\mathfrak{dg}_1)^{{}^{\vee}}\otimes gr_{1}^{\mathcal{D}}\mathbb{Q}\langle e^{0},e^{1},e^{-1}\rangle\\
\downarrow& &\downarrow\\
gr_{2}^{\mathcal{D}}\mathcal{H}&\stackrel{\partial}{\rightarrow}&gr_{1}^{\mathcal{D}}\mathcal{H}\otimes gr_{1}^{\mathcal{D}}\mathcal{H},
\end{array}$$
where
\[
\overline{\partial}=\sum_{n\geq 0} (\overline{\sigma_{2n+1}})^{{}^{\vee}}\otimes \overline{\partial}_{2n+1},
\]
and $(\overline{\sigma_{2n+1}})^{{}^{\vee}}, n\geq 0,$ is the dual basis of $\overline{\sigma_{2n+1}},n \geq 0,$ in $(\mathfrak{dg}_1)^{{}^{\vee}}$. The second column map transforms $(\overline{\sigma_{2n+1}})^{{}^{\vee}}\otimes \overline{\partial}_{2n+1}(x)$ to $\zeta^{\mathfrak{m}}(\overline{2n+1})\otimes \partial _{2n+1}(\mathcal{I}^{\mathfrak{m}}(x))$.

Thus we can calculate the image of $\mathcal{P}_{ev}^{o,\mathfrak{m}}$ under $\partial$ by calculating its lift on $$gr_{2}\mathbb{Q}\langle e^{0},e^{1},e^{-1}\rangle.$$

For even integers $n_{1}, n_{2}$, the motivic double zeta value $\zeta^{o,\mathfrak{m}}(n_{1},n_{2})$ regarded as an element of $gr_{2}\mathcal{H}$ is equal to
$$\frac{1}{4}\mathcal{I}^{\mathfrak{m}}(e^{1}(e^{0})^{n_{1}-1}e^{1}(e^{0})^{n_{2}-1}-e^{-1}(e^{0})^{n_{1}-1}e^{1}(e^{0})^{n_{2}-1}$$
$$-e^{-1}(e^{0})^{n_{1}-1}e^{-1}(e^{0})^{n_{2}-1}+e^{1}(e^{0})^{n_{1}-1}e^{-1}(e^{0})^{n_{2}-1}).$$
Denote the above expression by $\Lambda(n_{1}, n_{2})$. Let $\Theta_{i}^{n_{k}}$ be as above, $n>0$ and $s=n_{1}+n_{2}-2n-2$, according to Lemma \ref{4.1} and \ref{4.2} we have the following formula:

\ \ \ \ \ \ $\overline{\partial_{2n+1}}(\Lambda(n_{1}, n_{2}))=\frac{1}{4}(\Theta_{1}^{n_{1}}-\Theta_{1}^{n_{2}}-\Theta_{-1}^{n_{1}}+\Theta_{-1}^{n_{2}})$

\ \ \ \ \ \ \ \ \ \ \ \ \ \ \ \ \  \ \ \ \ \ \ \ \ \ \ \ $=\frac{1}{4}\left({2n \choose n_{1}-1}-{2n \choose n_{2}-1}\right)(e^{1}(e^{0})^{s}-e^{-1}(e^{0})^{s}).$

\

By Lemma \ref{2.1}, Lemma \ref{2.2}, if $s>0$, we have
$$\partial_{2n+1}(\zeta^{o,\mathfrak{m}}(n_{1},n_{2}))=\frac{1-2^{s+1}}{4(2^{s}-1)}\left({2n \choose n_{1}-1}-{2n \choose n_{2}-1}\right)\zeta^{\mathfrak{m}}(\overline{s+1}).$$

Combining with the definition of $\widetilde{D}$ and Proposition \ref{3.2} , it is obvious that $\mathrm{Im}\,(\widetilde{\partial})=\mathrm{Ker}\, (\widetilde{D})$.   $\hfill\Box$\\

Kaneko and Tasaka \cite{ref7} proved that there are at least $\mathrm{dim}\,\mathcal{S}_k(\Gamma_0(2))$-linear independent relations among the numbers $\{\zeta^{o}(k_1,k_2),k_1+k_2=k,k_1,k_2\geq 2, \mathrm{even}\}$.
From  Theorem \ref{ES} we obtain
\[
\mathrm{dim}_{\mathbb{Q}}\langle \zeta^{o}(k_1,k_2); k_1+k_2=k, k_1,k_2\geq 2, \mathrm{even}\rangle_{\mathbb{Q}}\leq \frac{k}{2}-1-\mathrm{dim}\,\mathcal{S}_k(\Gamma_0(2))
\]
immediately. The above inequality is compatible with Kaneko and Tasaka's result.

The next theorem gives an affirmative answer for part of  Kaneko and Tasaka's conjectures  in the motivic setting.

\begin{Thm}\label{5.2}
For an even integer $N\geq 4$. The elements
$$\{\zeta^{o,\mathfrak{m}} (k,N-k), 1\leq k\leq N-1\ odd\}$$
are $\mathbb{Q}$-linear independent. What's more, the set of their images in $gr_2^{\mathcal{D}}\mathcal{H}$  is a basis  of $gr_2^{\mathcal{D}}\mathcal{H}_N$.
\end{Thm}

\noindent{\bf Proof:}
We will make use of the above calculations again. The case $N=4$ is easy to check.
Given an even integer $N\geq 6$, according to Lemma \ref{4.3} and Lemma \ref{4.4}, for any odd $n_{1}$, $n_{2}$ such that $n_{1}+n_{2}=N$, we have for all $n_{1},n_2>1$, $\partial_{1}\left(\zeta^{\mathfrak{m}}\binom{n_{1},n_{2}}{\epsilon_1,\epsilon_2}\right)=0$ and
$$\partial_{1}(\zeta^{o,\mathfrak{m}}(1,n_{2}))=\frac{1}{2}\mathcal{I}^{\mathfrak{m}}[-e^{1}(e^{0})^{n_{2}-1}+e^{-1}(e^{0})^{n_{2}-1}],$$
$$\partial_{1}(\zeta^{o,\mathfrak{m}}(n_1,1))=\frac{1}{4}\mathcal{I}^{\mathfrak{m}}[e^{1}(e^{0})^{n_{1}-1}-e^{-1}(e^{0})^{n_{1}-1}].$$
Thus by the distribution formula, we have
$$\partial_{1}(\zeta^{o,\mathfrak{m}}(1,n_2))=\frac{1-2^{n_2}}{2-2^{n_{2}}}\zeta^{\mathfrak{m}}(\overline{n_{2}}),
$$
$$\partial_{1}(\zeta^{o,\mathfrak{m}}(n_1, 1)) =-\frac{1-2^{n_1}}{4(1-2^{n_1-1})}\zeta^{\mathfrak{m}}(\overline{n_1})                                                                 .$$

For the same reason, when $n\geq 1$, let $s=n_{1}+n_{2}-2n-2$, if $s>0$, the following formula holds:
\[
\begin{split}
&\;\;\;\;\partial_{2n+1}(\zeta^{o,\mathfrak{m}}(n_{1},n_{2}))\\
&=\frac{2^{s+1}-1}{4(2^s-1)}\left[{2n \choose n_{1}-1}-{2n \choose n_{2}-1}+\delta\binom{2n}{n_1-1}(1-2^{2n+1})\right]\zeta^{\mathfrak{m}}(\overline{s+1}).
\end{split}
\]
If $s=0$, we have
\[
\partial_{2n+1}(\zeta^{o,\mathfrak{m}}(n_{1},n_{2}))=-\frac{1}{4}\delta\binom{2n}{n_1-1}(2^{2n+1}-1)\zeta^{\mathfrak{m}}(\overline{1}).
\]

In conclusion, we can write the map $\partial$ in the following form in the case of $n_{1}+n_{2}=N$:
$$\partial\left(
            \begin{array}{c}
            \zeta^{o,\mathfrak{m}}(1,N-1) \\
            \zeta^{o,\mathfrak{m}}(3,N-3) \\
            \vdots                         \\
            \zeta^{o,\mathfrak{m}}(N-1,1) \\
            \end{array}
          \right)
=\widetilde{M} B \left(
    \begin{array}{c}
      \zeta^{\mathfrak{m}}(\overline{1})\otimes \zeta^{\mathfrak{m}}(\overline{N-1}) \\
      \zeta^{\mathfrak{m}}(\overline{3})\otimes \zeta^{\mathfrak{m}}(\overline{N-3}) \\
      \vdots \\
      \zeta^{\mathfrak{m}}(\overline{N-1})\otimes \zeta^{\mathfrak{m}}(\overline{1}) \\
    \end{array}
  \right).
$$

In the above formula $$B=\mathrm{diag}\left(\frac{1-2^{N-1}}{2-2^{N-1}}, \frac{2^{N-3}-1}{4(2^{N-4}-1)}, \cdots ,\frac{2^3-1}{4(2^2-1)},-\frac{(2^{N-1}-1)}{4}\right)$$ is a ${(\frac{N}{2})}^{th}$ invertible diagonal matrix, $\widetilde{M}$ is a square matrix of order $\frac{N}{2}$ in the form
\[
\widetilde{M}=\begin{pmatrix}
1 & &\cdots& &0  \\
0 & &{} && 0\\
\vdots && M && \vdots \\
0 & &{} && 0\\
c && \cdots && 1
\end{pmatrix},
\]
where $M$ is a $(\frac{N}{2}-2)^{th}$ square matrix in the middle of $\widetilde{M}$. The matrix $M=(a_{i,j})_{1\leq i,j \leq \frac{N}{2}-2}$,
\[
a_{i,j}=\binom{2j}{2i}-\binom{2j}{N-2-2i}+\delta\binom{2i}{2j}(1-2^{2j+1}).
\]

The theorem holds if ${M}$ is invertible by the fact that $\partial$ is injective.  $M$  can be written as the form $D+A$, where $t=\frac{N}{2}-1$,
$$D=\mathrm{diag}(d_{1},\cdots,d_{t-1}),\ A=(b_{i,j})_{1\leq i,j\leq t-1},$$
and
\[
d_i=1-2^{2i+1},\,\, b_{i,j}=\binom{2j}{2i}-\binom{2j}{2t-2i}.
\]

Given $j$, notice that $b_{i,j}+b_{t+1-i,j}=0$ and $b_{i,j}=0$ for $j<i<t-j$, it's obvious to check that
$$\sum_{ i=1}^{t-1}|b_{i,j}|=2\sum_{i=1}^{min\{\frac{t-1}{2},\ j-1\}}|b_{i,j}|\leq2\sum_{i=1}^{min\{\frac{t-1}{2},\ j-1\}}{2j \choose 2i}\leq 2\sum_{i=1}^{j-1}\binom{2j}{2i}< 2^{2j+1}-1.$$
So clearly for $j=1,...,t-1$, we have
\[
\sum_{i=1,i\neq j}^{t-1} | b_{i,j}|=\sum_{ i=1}^{t-1}|b_{i,j}|-|b_{j,j}|<|d_j+b_{j,j}|.
\]
By the following lemma, the matrix $M$, and furthermore $\widetilde{M}$ are invertible.$\hfill\Box$\\

\begin{lem}\label{dy}
For a real matrix $A=(a_{i,j})_{1\leq i,j\leq n}$, if $|a_{i,i}|> \sum_{i\neq j}|a_{i,j}|$ for $i=1,2,...,n$, then $|A|\neq 0$.
\end{lem}

\noindent{\bf Proof:}
Denote by $\alpha_{i}$ the $i^{th}$ column vector of $A$, if $|A|=0$, there exist $\{k_{1}, \cdots, k_{n}\}\neq \{0\}$ such that $k_{1}\alpha_{1}+\cdots+k_{n}\alpha_{n}=0$ is the zero column vector.

 Let
$$|k_{l}|=max \{|k_{1}|,\cdots,|k_{n}|\}.$$
Now consider the $l^{th}$  variable in the above zero column vector, because $|a_{l,l}|> \sum_{l\neq j}|a_{l,j}|$, $k_{1}a_{l,1}+\cdots+k_{n}a_{l,n}\neq 0$, we get a contradiction.
$\hfill\Box$\\

\begin{rem}
Kaneko and Tasaka \cite{ref7} conjectured that for given $N\geq 4$, elements $$\{\zeta^{o}(n_{1},n_{2});\ n_{1}\geq 1, n_{2}>1, odd,\ n_{1}+n_{2}=N\}$$ are $\mathbb{Q}$-linear independent. Theorem \ref{5.2} gives  a proof of the motivic version of Kaneko and Tasaka's conjecture.
\end{rem}

    As we know, for odd $n> 1$, the double zeta value $\zeta^{o}(n,1)$ is not well-defined. However, the motivic sum odd double zeta value $\zeta^{o,\mathfrak{m}}(n,1)$ is well-defined. We will calculate the period of $\zeta^{o,\mathfrak{m}}(n,1)$. Recall that
\[\tag{1}
\zeta^{o,\mathfrak{m}}(n,1)=\frac{1}{4}[\zeta^{\mathfrak{m}}(n,1)-\zeta^{\mathfrak{m}}(\overline{n},1)-\zeta^{\mathfrak{m}}(n,\overline{1})+\zeta^{\mathfrak{m}}(\overline{n},\overline{1})].
\]

\begin{lem}\label{per}
For $n>1,\mathrm{odd}$, the period of $\zeta^{o,\mathfrak{m}}(n,1)$ is
\[
\begin{split}
&\,\,\,\,\,\,\,\,per(\zeta^{o,\mathfrak{m}}(n,1))\\
&=\frac{1}{4}[-\zeta(1,n)-\zeta(n+1)-\zeta(\overline{1},n)+\zeta(\overline{n+1})-\zeta(n,\overline{1})+\zeta(\overline{n},\overline{1})]\\
&=\frac{1}{4}[  -\zeta(1,n)+\zeta(1, \overline{n})-\zeta(n,\overline{1})+\zeta(\overline{n},\overline{1}) +(2^{-n}-2)\zeta(n+1)               ].
\end{split}
\]
\end{lem}

\noindent{\bf Proof:}
It is direct to get the periods of $\zeta^{\mathfrak{m}}(n,\overline{1})$ and $\zeta^{\mathfrak{m}}(\overline{n},\overline{1})$, we only need to determine the other two terms in formula $(1)$. Consider the following regularized integral:
$$\int_{0<t_{1}<\cdots <t_{n+1}<1-\eta}\frac{dt_{1}}{1-t_{1}}\frac{dt_{2}}{t_{2}}\cdots \frac{dt_{n}}{t_{n}}\frac{dt_{n+1}}{1-t_{n+1}}$$
$$=\sum_{0<s<r}\frac{(1-\eta)^{r}}{s^{n}r}=\sum_{s=1}^{\infty}\frac{1}{s^{n}}\left(-log(\eta)-\sum_{r=1}^{s}\frac{(1-\eta)^{r}}{r}\right).$$

Let $log(\eta)=0$, the above integral equals to $-\sum_{s=1}^{\infty}\sum_{r=1}^{s}\frac{(1-\eta)^{r}}{rs^{n}}$, and then let $\eta\rightarrow 0$ we have
$$-\sum_{s=1}^{\infty}\sum_{r=1}^{s}\frac{1}{rs^{n}}=-\sum_{0<r<s}\frac{1}{rs^{n}}-\sum_{s=1}^{\infty}\frac{1}{s^{n+1}}=-\zeta(1,n)-\zeta(n+1).$$
By the definition of $per$, we have
\[
per(\zeta^{\mathfrak{m}}(n,1))=-\zeta(1,n)-\zeta(n+1).
\]
The same as the above calculation, we have
$$per(\zeta^{\mathfrak{m}}(\overline{n},1))=-\zeta(1,\overline{n})-\zeta(\overline{n+1}).$$
Combining with the formula $(1)$, lemma proved.$\hfill\Box$\\

The following remark follows from Theorem \ref{5.2}  and Lemma  \ref{per} immediately.

\begin{rem}
Every element $\zeta\binom{n_1,n_2}{\epsilon_1,\epsilon_2}, n_1+n_2=N,N\, \mathrm{even}, (n_2,\epsilon_2)\neq (1,1)$ can be written as a  $\mathbb{Q}$-linear combination of $\zeta^{o}(odd,odd),\zeta(N)$ and $\mathrm{per}(\zeta^{o,\mathfrak{m}}(N-1,1))$ as above.
\end{rem}

\subsection{Kaneko and Tasaka's three conjectures}
 Kaneko and Tasaka \cite{ref7} also conjectured that $\langle\zeta^{o}(r,N-r); 1\leq r\leq N-2\rangle_{\mathbb{Q}}$ is spanned by $\zeta^{o}(odd,odd)$ and $\zeta(N)$. We will  prove this statement as an application of the motivic method.

\begin{Thm}\label{span}
For a given even integer $N\geq 4$, the space $\langle\zeta^{o}(r, N-r); 1\leq r\leq N-2\rangle_{\mathbb{Q}}$ is spanned by
$$\{\zeta(N),\zeta^{o}(r, N-r); 1\leq r\leq N-3,\mathrm{odd}\}.$$
\end{Thm}

\noindent{\bf Proof:}
Denote by $gr_2^{\mathcal{D}}\mathcal{H}_N$ the weight $N$ part of $gr_2^{\mathcal{D}}\mathcal{H}$.
According to the property of the period map $per$, we only need to prove that
$$\langle\zeta^{o,\mathfrak{m}}(r, N-r); 1\leq r\leq N-2\rangle_{\mathbb{Q}}=\mathrm{span}\{\zeta^{o,\mathfrak{m}}(r, N-r); 1\leq r\leq N-3,\mathrm{ odd}\}$$
in $gr_2^{\mathcal{D}}\mathcal{H}_N$. (Be ware that $\zeta^{o,\mathfrak{m}}(r,N-r);1\leq r\leq N-3,\mathrm{ odd}$ are elements of $\mathcal{H}$, in the above formula we  mean their natural images in $gr_2^{\mathcal{D}}\mathcal{H}_N$.)

We use the same notation as in the proof of Theorem \ref{5.2}, there is a matrix $E$ such that
$$\partial\left(
            \begin{array}{c}
            \zeta^{o,\mathfrak{m}}(1,N-1) \\
            \zeta^{o,\mathfrak{m}}(3,N-3) \\
            \vdots                         \\
            \zeta^{o,\mathfrak{m}}(N-1,1) \\
            \end{array}
          \right)
=E\left(
    \begin{array}{c}
      \zeta^{\mathfrak{m}}(\overline{1})\otimes \zeta^{\mathfrak{m}}(\overline{N-1}) \\
      \zeta^{\mathfrak{m}}(\overline{3})\otimes \zeta^{\mathfrak{m}}(\overline{N-3}) \\
      \vdots \\
      \zeta^{\mathfrak{m}}(\overline{N-1})\otimes \zeta^{\mathfrak{m}}(\overline{1}) \\
    \end{array}
  \right),
$$
where $E=\widetilde{M} B$ is invertible.

Thus  we have
$$\partial E^{-1}\left(
            \begin{array}{c}
            \zeta^{o,\mathfrak{m}}(1,N-1) \\
            \zeta^{o,\mathfrak{m}}(3,N-3) \\
            \vdots                         \\
            \zeta^{o,\mathfrak{m}}(N-1,1) \\
            \end{array}
          \right)
=\left(
    \begin{array}{c}
      \zeta^{\mathfrak{m}}(\overline{1})\otimes \zeta^{\mathfrak{m}}(\overline{N-1}) \\
      \zeta^{\mathfrak{m}}(\overline{3})\otimes \zeta^{\mathfrak{m}}(\overline{N-3}) \\
      \vdots \\
      \zeta^{\mathfrak{m}}(\overline{N-1})\otimes \zeta^{\mathfrak{m}}(\overline{1}) \\
    \end{array}
  \right).
$$

On the other hand, according to Lemma \ref{4.1} and Lemma \ref{4.2}, we have
$$\partial\left(
            \begin{array}{c}
            \zeta^{o,\mathfrak{m}}(2,N-2) \\
            \zeta^{o,\mathfrak{m}}(4,N-4) \\
            \vdots                         \\
            \zeta^{o,\mathfrak{m}}(N-2,2) \\
            \end{array}
          \right)
=F\left(
    \begin{array}{c}
      \zeta^{\mathfrak{m}}(\overline{1})\otimes \zeta^{\mathfrak{m}}(\overline{N-1}) \\
      \zeta^{\mathfrak{m}}(\overline{3})\otimes \zeta^{\mathfrak{m}}(\overline{N-3}) \\
      \vdots \\
      \zeta^{\mathfrak{m}}(\overline{N-1})\otimes \zeta^{\mathfrak{m}}(\overline{1}) \\
    \end{array}
  \right),
$$
where $F$ is a matrix of order $(\frac{N}{2}-1,\frac{N}{2})$, thus
$$\partial\left(
            \begin{array}{c}
            \zeta^{o,\mathfrak{m}}(2,N-2) \\
            \zeta^{o,\mathfrak{m}}(4,N-4) \\
            \vdots                         \\
            \zeta^{o,\mathfrak{m}}(N-2,2) \\
            \end{array}
          \right)
=\partial FE^{-1}\left(
            \begin{array}{c}
            \zeta^{o,\mathfrak{m}}(1,N-1) \\
            \zeta^{o,\mathfrak{m}}(3,N-3) \\
            \vdots                         \\
            \zeta^{o,\mathfrak{m}}(N-1,1) \\
            \end{array}
          \right).
$$

By the injectivity of $\partial$, we have
$$\left(
            \begin{array}{c}
            \zeta^{o,\mathfrak{m}}(2,N-2) \\
            \zeta^{o,\mathfrak{m}}(4,N-4) \\
            \vdots                         \\
            \zeta^{o,\mathfrak{m}}(N-2,2) \\
            \end{array}
          \right)
= FE^{-1}\left(
            \begin{array}{c}
            \zeta^{o,\mathfrak{m}}(1,N-1) \\
            \zeta^{o,\mathfrak{m}}(3,N-3) \\
            \vdots                         \\
            \zeta^{o,\mathfrak{m}}(N-1,1) \\
            \end{array}
          \right)
$$
in $gr_2^{\mathcal{D}}\mathcal{H}_N$.
From the explicit calculation in Theorem \ref{ES} and Theorem \ref{5.2}, it's obvious to check that the last column of the matrix $FE^{-1}$ is $0$. By using the period map, the theorem is proved.
$\hfill\Box$\\

Kaneko and Tasaka \cite{ref7} gave some other conjectures and we can prove them by the same motivic method as above.

\begin{prop}
For even integer $N\geq6$, we have
$$\langle\zeta(n_{1},n_{2}); n_{1}+n_{2}=N, n_{2}\geq 2\rangle_{\mathbb{Q}}\subseteq \langle\zeta^{o}(n_{1},n_{2}); n_{1}+n_{2}=N, n_{i} \geq 2\rangle_{\mathbb{Q}},$$
	$$\langle\zeta^{o}(n_{1},n_{2}); n_{1}+n_{2}=N, n_{i}\ even\rangle_{\mathbb{Q}}\subseteq \langle\zeta(n_{1},n_{2}); n_{1}+n_{2}=N, n_{2}\geq 2\rangle_{\mathbb{Q}}.$$
\end{prop}
\noindent{\bf Proof:}
We only need to prove this proposition in the motivic version. According to our calculations above, for $n\geq 0$, let $s=N-2n-2$, we have
$$\partial_{1}(\zeta^{\mathfrak{m}}(n_{1},n_{2}))=0,$$
\[
\begin{split}
&\;\;\;\;\partial_{2n+1}(\zeta^{\mathfrak{m}}(n_{1},n_{2}))\\
&=2^{2n}\left((-1)^{n_{1}}{2n \choose n_{1}-1}-(-1)^{n_{2}}{2n \choose n_{2}-1}+\delta\binom{2n}{n_{1}-1}\right)\zeta^{\mathfrak{m}}(s+1).
\end{split}
\]
By the distribution formula, when $s\neq 0$ we have
\[
\begin{split}
&\;\;\;\;\partial_{2n+1}(\zeta^{\mathfrak{m}}(n_{1},n_{2}))\\
&=\frac{2^{N-2}}{1-2^{s}}\left((-1)^{n_{1}}{2n \choose n_{1}-1}-(-1)^{n_{2}}{2n \choose n_{2}-1}+\delta\binom{2n}{n_{1}-1}\right)\zeta^{\mathfrak{m}}(\overline{s+1})
\end{split}
\]
and when $s=0$ we have
$$\partial_{2n+1}(\zeta^{\mathfrak{m}}(n_{1},n_{2}))=0.$$

 We have shown there is an invertible matrix $S$ such that
$$\partial S^{-1}\left(
            \begin{array}{c}
            \zeta^{o,\mathfrak{m}}(3,N-3) \\
            \vdots                         \\
            \zeta^{o,\mathfrak{m}}(N-3,3) \\
            \end{array}
          \right)
=\left(
    \begin{array}{c}
      \zeta^{\mathfrak{m}}(\overline{3})\otimes \zeta^{\mathfrak{m}}(\overline{N-3}) \\
      \vdots \\
      \zeta^{\mathfrak{m}}(\overline{N-3})\otimes \zeta^{\mathfrak{m}}(\overline{3}) \\
    \end{array}
  \right).
$$
By the injectivity of $\partial$ , we have
$$\langle\zeta^{\mathfrak{m}}(n_{1},n_{2}); n_{1}+n_{2}=N, n_{2}\geq 2\rangle_{\mathbb{Q}}\subseteq \langle\zeta^{o,\mathfrak{m}}(n_{1},n_{2}); n_{1}+n_{2}=N, 2\leq n_{i}\leq k-2\rangle_{\mathbb{Q}}.$$

For $n_1,n_2\geq 2$, even,  if $s=0$ or $n=0$,
\[
\partial_{2n+1}(\zeta^{o,\mathfrak{m}}(n_1,n_2)   )=0,
\]
 and if $s,n>0$,
\[
\begin{split}
\;\;\;\;\partial_{2n+1}(\zeta^{o,\mathfrak{m}}(n_1,n_2))
=\frac{1}{4}\frac{1-2^{s+1}}{2^s-1}\left[\binom{2n}{n_1-1}-\binom{2n}{n_2-1}   \right]\zeta^{\mathfrak{m}}(\overline{s+1} ).
\end{split}
\]
Since the map $\partial$ is injective, to prove
$$\langle\zeta^{o}(n_{1},n_{2}); n_{1}+n_{2}=N, n_{i}\ even\rangle_{\mathbb{Q}}\subseteq \langle\zeta(n_{1},n_{2}); n_{1}+n_{2}=N, n_{2}\geq 2\rangle_{\mathbb{Q}},$$
it suffices to prove that there are numbers $d\binom{m_1,m_2}{n_1,n_2}$ which satisfy
\[
\begin{split}
&\;\;\;\;\frac{1}{2^{2n}}\left[\binom{2n}{n_1-1}-\binom{2n}{n_2-1}\right]\\
&=\sum_{\substack{m_1+m_2=N\\m_i\geq 1}}d\binom{m_1,m_2}{n_1,n_2}
\left[ (-1)^{m_1}\binom{2n}{m_1-1}-(-1)^{m_2}\binom{2n}{m_2-1}+\delta\binom{2n}{m_1-1} \right]
\end{split}
\]
for all $n_1+n_2=N$, $n_1,n_2\geq 2$, even, $3\leq 2n+1 \leq N-3$.
The above statement follows from Lemma \ref{inver} and Remark \ref{calinv} below.
$\hfill\Box$\\

Denote by
\[
\begin{split}
&V_{N,2}=\langle x_1^{n_1-1}x_2^{n_2-1};n_1+n_2=N, n_1,n_2\geq 3,\mathrm{odd}\rangle_{\mathbb{Q}},\\
&\mathbb{P}_{N,2}=\langle x_1^{n_1-1}x_2^{n_2-1};n_1+n_2=N, n_1,n_2\geq 1\rangle_{\mathbb{Q}},\\
&\mathbb{P}_{N,2}^{od}=\langle x_1^{n_1-1}x_2^{n_2-1};n_1+n_2=N, n_1,n_2\geq 2,\mathrm{even}\rangle_{\mathbb{Q}}.
\end{split}
\]
For $p(x_1,x_2)\in V_{N,2}$, define $$L_{1,1}(p)(x_1,x_2)=p(x_1,x_2)+p(x_1-x_2,x_1)-p(x_1-x_2,x_2),$$
$$L_{\frac{1}{2},1}(p)(x_1,x_2)=p(\frac{x_1}{2},x_2)+p(\frac{x_1-x_2}{2},x_1)-p(\frac{x_1-x_2}{2},x_2).$$
Then we have:
\begin{lem}\label{inver}
Denote by $i^{od}: \mathbb{P}_{N,2}\to \mathbb{P}^{od}_{N,2}$ the natural map which satisfies for $p(x_1,x_2)\in \mathbb{P}_{N,2}$, $$i^{od}(p)(x_1,x_2)=p(x_1,x_2)-p(-x_1,x_2).$$
There is a linear map $j: \mathbb{P}_{N,2}\to \mathbb{P}^{od}_{N,2}$ such that the following diagram is commutative
\[
 \xymatrix{
   V_{N,2} \ar[d]^{L_{\frac{1}{2},1}} \ar[r]^{ L_{1,1}} & \mathbb{P}_{N,2} \ar[d]^{j} \\
 \mathbb{P}_{N,2} \ar[r]^{i^{od}} & \mathbb{P}_{N,2}^{od} . }
 \]
\end{lem}
\noindent{\bf Proof:}
Define $j_1: \mathbb{P}_{N,2}\to \mathbb{P}_{N,2}$ as the $\mathbb{Q}$-linear map which is induced by
\[
x_1\mapsto \frac{x_1+x_2}{2}, x_2\mapsto x_2.
\]
Define $j_2: \mathbb{P}_{N,2}\to \mathbb{P}_{N,2}$ as the $\mathbb{Q}$-linear map which is induced by
\[
x_1\mapsto \frac{x_1+x_2}{2}, x_2\mapsto x_1.
\]
Define $j=\frac{1}{2}i^{od}\circ(j_1-j_2)$.

For $p\in V_{N,2}$,
\[
\begin{split}
&\;\;\;\;\;i^{od}\circ L_{\frac{1}{2},1}(p)(x_1,x_2)\\
&= L_{\frac{1}{2},1}(p)(x_1,x_2)- L_{\frac{1}{2},1}(p)(-x_1,x_2)\\
&=p(\frac{x_1-x_2}{2},x_1)-p(\frac{x_1-x_2}{2},x_2)-p(\frac{x_1+x_2}{2},x_1)+p(\frac{x_1+x_2}{2},x_2),
\end{split}
\]
\[
\begin{split}
&\;\;\;\;\;j\circ L_{1,1}(p)(x_1,x_2)\\
&=\frac{1}{2}\left[(j_1-j_2)\circ L_{1,1}(p)(x_1,x_2)-  (j_1-j_2)\circ L_{1,1}(p)(-x_1,x_2)\right] \\
& =\frac{1}{2}(j_1\circ L_{1,1}(p)(x_1,x_2) -j_2\circ L_{1,1}(p)(x_1,x_2)                                                  -  j_1\circ L_{1,1}(p)(-x_1,x_2)\\
&\;\;\;\;\;\;\;+ j_2\circ L_{1,1}(p)(-x_1,x_2) )\\
&=\frac{1}{2}(L_{1,1}(p)(\frac{x_1+x_2}{2},x_2)-L_{1,1}(p)(\frac{x_1+x_2}{2},x_1)-
L_{1,1}(p)(\frac{-x_1+x_2}{2},x_2)\\
&\;\;\;\;\;\;\;+ L_{1,1}(p)(\frac{-x_1+x_2}{2},-x_1)) \\
&=\frac{1}{2}\left[p(\frac{x_1+x_2}{2},x_2)+p(\frac{x_1-x_2}{2},\frac{x_1+x_2}{2})-p(\frac{x_1-x_2}{2},x_2)\right]\\
&\;\;\;\;-\frac{1}{2}\left[p(\frac{x_1+x_2}{2},x_1)+p(\frac{x_1-x_2}{2},\frac{x_1+x_2}{2})-p(\frac{x_1-x_2}{2},x_1)\right]\\
&\;\;\;\;-\frac{1}{2}\left[p(\frac{-x_1+x_2}{2},x_2)+p(\frac{x_1+x_2}{2},\frac{x_1-x_2}{2})-p(\frac{x_1+x_2}{2},x_2)\right]\\
&\;\;\;\;+\frac{1}{2}\left[p(\frac{-x_1+x_2}{2},x_1)+p(\frac{x_1+x_2}{2},\frac{x_1-x_2}{2})-p(\frac{x_1+x_2}{2},x_1)\right]\\
&=p(\frac{x_1-x_2}{2},x_1)-p(\frac{x_1-x_2}{2},x_2)-p(\frac{x_1+x_2}{2},x_1)+p(\frac{x_1+x_2}{2},x_2).
\end{split}
\]
As a result of the above calculations, the lemma is proved.
$\hfill\Box$\\
\begin{rem}\label{calinv}
Define $d\binom{m_1,m_2}{n_1,n_2}$ as the coefficient of $x_1^{n_1-1}x_2^{n_2-1}$ in $$\frac{1}{2}j(x_1^{m_1-1}x_2^{m_2-1}),$$
i.e.
\[
\frac{1}{2}j(x_1^{m_1-1}x_2^{m_2-1})=\sum_{\substack{n_1+n_2=N\\n_i\geq 2, even}}d\binom{m_1,m_2}{n_1,n_2} x_1^{n_1-1}x_2^{n_2-1}.
\]
For $3\leq 2n+1\leq N-3$, by running the commutative diagram in Lemma \ref{inver} on $$p=x_1^{2n}x_2^{N-2-2n}\in V_{N,2},$$
we have
  \[
\begin{split}
&\;\;\;\;\frac{1}{2^{2n}}\left[\binom{2n}{n_1-1}-\binom{2n}{n_2-1}\right]\\
&=\sum_{\substack{m_1+m_2=N\\m_i\geq 1}}d\binom{m_1,m_2}{n_1,n_2}
\left[ (-1)^{m_1}\binom{2n}{m_1-1}-(-1)^{m_2}\binom{2n}{m_2-1}+\delta\binom{2n}{m_1-1} \right].
\end{split}
\]
\end{rem}

  By the motivic method we can also prove the following proposition which was proved by Kaneko and Tasaka \cite{ref7}.

\begin{prop}
For odd integer $N>6$, we have
$$\langle\zeta^{o}(n_{1},n_{2}); n_{1}+n_{2}=N, n_{i}\geq 2\rangle_{\mathbb{Q}}\subseteq \langle\zeta(n_{1},n_{2}); n_{1}+n_{2}=N, n_{2}\geq 2\rangle_{\mathbb{Q}}.$$
\end{prop}

\section{The higher depth case}\label{hdc}

 In this section we calculate the map $$\partial:gr_r^{\mathcal{D}}\mathcal{H} \rightarrow gr_{1}^{\mathcal{D}}\mathcal{H}^{odd}\otimes gr_{r-1}^{\mathcal{D}}\mathcal{H}$$ in the case of $r\geq 3$ for sum odd motivic multiple zeta values  explicitly.  As a corollary we obtain a basis for the depth-graded motivic triple zeta values of odd weight.  What's more, all elements of this basis are the natural images of sum odd motivic multiple zeta values in the depth-graded motivic triple zeta values of odd weight. At last we conjecture that a matrix appeared in the explicit calculation of $\partial$ on the sum odd motivic multiple zeta values is invertible.

 Denote by
$$T_{N,r}=\{(n_{1},\cdots,n_{r})\in \mathbb{Z}^{r};n_{1}+\cdots+n_{r}=N, n_{i}\geq1,\mathrm{odd},1\leq i\leq r\}.$$
Define $\delta\binom{m_1,\cdots,m_r}{n_1,\cdots,n_r}=1$ if $(m_1,\cdots,m_r)=(n_1,\cdots,n_r)$, $\delta\binom{m_1,\cdots,m_r}{n_1,\cdots,n_r}=0$ if $(m_1,\cdots,m_r)\neq (n_1,\cdots,n_r)$.
\begin{prop}\label{tran}
Let $N\equiv r\ mod\ 2$ and $N\geq r+2$, for $(k_{1},\cdots,k_{r})\in T_{N,r}$ we have
$$\partial(\zeta^{o,\mathfrak{m}}(k_{1},\cdots,k_{r}))=\sum_{(n_{1},\cdots,n_{r})\in T_{N,r}}e\binom{k_1,k_2,\cdots,k_r}{n_1,n_2,\cdots,n_r}\zeta^{\mathfrak{m}}(\overline{n_{1}})\otimes\zeta^{o,\mathfrak{m}}(n_{2},\cdots,n_{r}),
$$
where for $n_{1}\geq 3,\mathrm{odd}$,
\[
\begin{split}
&e\binom{k_1,k_2,\cdots,k_r}{n_1,n_2,\cdots,n_r}=(2^{n_{1}-1}-\frac{1}{2})\delta\binom{k_1,k_2,\cdots,k_r}{n_1,n_2,\cdots,n_r}            \\
                                               &\;\;\;\;\;\;\;\;\;\;\;\;\;\;\;\;\;\;+\frac{1}{2}\sum_{i=1}^{r-1}\left(\binom{n_1-1}{k_{i+1}-1}-\binom{n_1-1}{k_i-1} \right)\delta\binom{k_1,\cdots,k_{i-1},k_{i+2},\cdots,k_r}{n_2,\cdots,n_i,\,\,n_{i+2},\cdots,n_r}
\end{split}
\]
and
$$e\binom{k_1,k_2,\cdots,k_r}{1,n_2,\cdots,n_r}=-\delta\binom{k_1, k_2,\cdots,k_r }{1, n_2,\cdots,n_r }+\frac{1}{2}\delta\binom{k_1,\cdots,k_{r-1},k_r }{n_2,\cdots,n_r, 1 }. $$
\end{prop}
\noindent{\bf Proof:}
Notice the following calculation:
\[
\begin{split}
&\,\,\,\,\zeta^{o,\mathfrak{m}}(n_{1},\cdots,n_{r})\\
&=\frac{1}{2^{r}}\sum_{\epsilon_i\in \{\pm1\},1\leq i\leq r}\epsilon_1\cdots \epsilon_r\zeta^{\mathfrak{m}}\binom{n_1,\cdots,n_r}{\epsilon_1,\cdots, \epsilon_r}\\
                                                                   &=\frac{1}{2^r}\sum_{\epsilon_i\in \{\pm1\},1\leq i\leq r}\epsilon_1\cdots \epsilon_r\mathcal{I}^{\mathfrak{m}}(e^{\epsilon_1\cdots \epsilon_r}(e^0)^{n_1-1}e^{\epsilon_2\cdots \epsilon_r}(e^0)^{n_2-1}\cdots e^{\epsilon_r}(e^0)^{n_r-1})\\
                                                                   &=\frac{1}{2^r}\sum_{\epsilon_i\in \{\pm 1\},2\leq i\leq r}\mathcal{I}^{\mathfrak{m}}[(e^{1}(e^0)^{n_1-1}e^{\epsilon_2\cdots \epsilon_r}(e^0)^{n_2-1}\cdots e^{\epsilon_r}(e^0)^{n_r-1}) -\\
                                                                   &\,\,\,\,\,\,\,\,\,\,\,\,\,\,\,\,\,\,\,\,\,\,\,\,\,\,\,\,\,\,\,\,\,\,\,\,\,\,\,\,\,\,\,\,e^{-1}(e^0)^{n_1-1}e^{\epsilon_2\cdots \epsilon_r}(e^0)^{n_2-1}\cdots e^{\epsilon_r}(e^0)^{n_r-1}) ]\\
                                                                   &=\frac{1}{2^r}\sum_{\epsilon_i\in \{\pm 1\},2\leq i\leq r}\mathcal{I}^{\mathfrak{m}}[(e^{1}(e^0)^{n_1-1}e^{\epsilon_2}(e^0)^{n_2-1}\cdots e^{\epsilon_r}(e^0)^{n_r-1}) -\\
                                                                   &\,\,\,\,\,\,\,\,\,\,\,\,\,\,\,\,\,\,\,\,\,\,\,\,\,\,\,\,\,\,\,\,\,\,\,\,\,\,\,\,\,\,\,\,e^{-1}(e^0)^{n_1-1}e^{\epsilon_2}(e^0)^{n_2-1}\cdots e^{\epsilon_r}(e^0)^{n_r-1}) ].
\end{split}
\]
 We have
$$\overline{\sigma_{2n+1}}\circ(e_{0}^{a_{0}}e_{i_{1}}e_{0}^{a_{1}}\cdots e_{i_{r}}e_{0}^{a_{r}})=
\sum_{j=1}^{r}e_{0}^{a_{0}}\cdots(\overline{\sigma_{2n+1}}\circ e_{i_{j}})e_{0}^{a_{j}}\cdots e_{i_{r}}e_{0}^{a_{r}}+ e_{0}^{a_{0}}
\cdots e_{i_{r}}e_{0}^{a_{r}}\overline{\sigma_{2n+1}}.$$

Since when $n=0$, we have
\[
\overline{\sigma}_1=e_{-1},e_{-1}\circ e_1=e_{-1}e_1-e_1e_{-1}, e_{-1}\circ e_{-1}=e_1e_{-1}-e_{-1}e_1.
\]
It follows that
\[
\begin{split}
&\,\,\,\,\overline{\partial}_1\left(e^{i_1}(e^0)^{a_1}\cdots e^{i_s}(e^0)^{a_s}\right)\\
&=\delta\binom{a_1}{0}\delta\binom{i_{1}i_2}{-1}i_1(e^{-1}-e^1)(e^0)^{a_2}\cdots e^{i_s}(e^0)^{a_s}\\
&\,\,\,\,\,\,\,+\cdots \\
&\,\,\,\,\,\,\,+\delta\binom{a_{s-1}}{0}\delta\binom{i_{s-1}i_s}{-1}i_{s-1}e^{i_1}(e^0)^{a_1}\cdots e^{i_{s-2}}(e^0)^{a_{s-2}}(e^{-1}-e^1)(e^0)^{a_s}\\
&\,\,\,\,\,\,\,+ \delta\binom{a_s}{0}\delta\binom{i_s}{-1}e^{i_1}(e^0)^{a_1}\cdots e^{i_{s-1}}(e^0)^{a_{s-1}}.
\end{split}
\]
As a result, we have
\[
\begin{split}
&\;\;\;\;\partial_1\left(\zeta^{o,\mathfrak{m}}(n_1,n_2,\cdots,n_r)\right)\\
&=\frac{1}{2^r}\sum_{\epsilon_i\in \{\pm 1\}, 1\leq i\leq r}\mathcal{I}^{\mathfrak{m}}[\overline{\partial}_1
\left(\epsilon_1e^{\epsilon_1}(e^0)^{n_1-1}e^{\epsilon_2}(e^0)^{n_2-1}\cdots e^{\epsilon_r}(e^0)^{n_r-1} \right)]\\
&=\frac{1}{2^r}  \delta\binom{n_1}{1}\sum_{\epsilon_i\in \{\pm 1\}, 1\leq i\leq r}\mathcal{I}^{\mathfrak{m}}
\left( \delta\binom{\epsilon_1\epsilon_2}{-1}\epsilon_1^2(e^{-1}-e^1)(e^0)^{n_2-1}\cdots e^{\epsilon_r}(e^0)^{n_r-1} \right)\\
&+\cdots\\
&+\frac{1}{2^r}  \delta\binom{n_{r-1}}{1}\sum_{\epsilon_i\in \{\pm 1\}, 1\leq i\leq r}\mathcal{I}^{\mathfrak{m}}[ \delta\binom{\epsilon_{r-1}\epsilon_r}{-1} \epsilon_1\epsilon_{r-1}e^{\epsilon_1}(e^0)^{n_1-1}  \\
&\,\,\,\,\,\,\,\,\,\,\,\,\,\,\,\,\,\,\,\,\,\,\,\,\,\,\,\,\,\,\,\,\,\,\,\,\,\,\,\,\,\,\,\,\,\,\,\,\,\,\,\,\,\,\,\,\,\,\,\,\,\,\,\,\,\,\,\,\,\,\,\,\,\,\,\,\,\,\, \cdots e^{\epsilon_{r-2}}(e^0)^{n_{r-2}-1}(e^{-1}-e^1)(e^0)^{n_r-1} ]\\
&+\frac{1}{2^r} \delta\binom{n_{r}}{1}\sum_{\epsilon_i\in \{\pm 1\}, 1\leq i\leq r}\mathcal{I}^{\mathfrak{m}}
\left(\delta\binom{\epsilon_r}{-1}\epsilon_1e^{\epsilon_1}(e^0)^{n_1-1}\cdots e^{\epsilon_{r-1}}(e^0)^{n_{r-1}-1} \right)\\
&=-\delta\binom{n_1}{1}\zeta^{o,\mathfrak{m}}(n_2,\cdots,n_r)
+\frac{1}{2}\delta\binom{n_r}{1}\zeta^{o,\mathfrak{m}}(n_1,\cdots,n_{r-1}).
\end{split}
\]

In the above calculation, the last equality is due to the fact that
\[
\sum_{\epsilon_1,\epsilon_2\in\{\pm 1\}}\delta\binom{\epsilon_1\epsilon_2}{-1}\epsilon_1=0.
\]
Similarly for $n\geq 1$, from
\[
\overline{\sigma_{2n+1}}=(1-2^{2n})\sum_{r=0}^{2n}(-1)^r{2n \choose r}e_{0}^{2n-r}e_{-1}e_{0}^{r}+2^{2n}\sum_{r=0}^{2n}(-1)^r{2n \choose r}e_{0}^{2n-r}e_{1}e_{0}^{r},\]

\[
\begin{split}
 \overline{\sigma_{2n+1}}\circ e_1=&(1-2^{2n})    \sum_{r=0}^{2n}(-1)^r{2n \choose r}\left(e_0^{2n-r}e_{-1}e_0^re_1-e_1e_0^{2n-r}e_{-1}e_0^r              \right)  \\
 &+2^{2n}   \sum_{r=0}^{2n}(-1)^r{2n \choose r}\left(   e_0^{2n-r}e_1e_0^re_1-e_1e_0^{2n-r}e_1e_0^r         \right),
 \end{split}
\]

\[
\begin{split}
 \overline{\sigma_{2n+1}}\circ e_{-1}=&(1-2^{2n})    \sum_{r=0}^{2n}(-1)^r{2n \choose r}\left(e_0^{2n-r}e_{1}e_0^re_{-1}-e_{-1}e_0^{2n-r}e_{1}e_0^r              \right)  \\
 &+2^{2n}   \sum_{r=0}^{2n}(-1)^r{2n \choose r}\left(   e_0^{2n-r}e_{-1}e_0^re_{-1}-e_{-1}e_0^{2n-r}e_{-1}e_0^r         \right),
 \end{split}
\]
we have
\[
\begin{split}
&\;\;\;\;\overline{\partial}_{2n+1}\left(e^{i_1}(e^0)^{a_1}\cdots e^{i_s}(e^0)^{a_s}\right)\\
&=(1-2^{2n})\delta\binom{a_1}{2n}\delta\binom{i_1i_2}{-1}e^{-i_1}(e^0)^{a_2}e^{i_3}(e^0)^{a_3}\cdots e^{i_s}(e^0)^{a_s}\\
\end{split}
\]
\[
\begin{split}
&-(1-2^{2n})(-1)^{a_1}\binom{2n}{a_1}\delta\binom{i_1i_2}{-1}e^{i_1}(e^0)^{a_1+a_2-2n}e^{i_3}(e^0)^{a_3}\cdots e^{i_s}(e^0)^{a_s}\\
&+2^{2n}\delta\binom{a_1}{2n}\delta\binom{i_1i_2}{1}e^{i_1}(e^0)^{a_2}e^{i_3}(e^{0})^{a_3}\cdots e^{i_s}(e^0)^{a_s}\\
&-2^{2n}(-1)^{a_1}\binom{2n}{a_1}\delta\binom{i_1i_2}{1}e^{i_1}(e^0)^{a_1+a_2-2n}e^{i_3}(e^0)^{a_3}\cdots e^{i_s}(e^0)^{a_s}\\
&+\cdots\\
&+(1-2^{2n})(-1)^{a_{s-1}}\binom{2n}{a_{s-1}}\delta\binom{i_{s-1}i_s}{-1}e^{i_1}(e^0)^{a_1}\cdots e^{i_{s-2}}(e^0)^{a_{s-2}+a_{s-1}-2n}e^{-i_{s-1}}(e^0)^{a_s} \\
&-(1-2^{2n})(-1)^{a_{s-1}}\binom{2n}{a_{s-1}}\delta\binom{i_{s-1}i_s}{-1} e^{i_1}(e^0)^{a_1}\cdots e^{i_{s-2}}            (e^0)^{a_{s-2}}e^{i_{s-1}}(e^0)^{a_{s-1}+a_s-2n} \\
&+2^{2n}(-1)^{a_{s-1}}\binom{2n}{a_{s-1}}\delta\binom{i_{s-1}i_s}{1}e^{i_1}(e^0)^{a_1}\cdots e^{i_{s-2}}(e^0)^{a_{s-2}+a_{s-1}-2n}e^{i_{s-1}}(e^0)^{a_s}\\
&+2^{2n}(-1)^{a_{s-1}}\binom{2n}{a_{s-1}}\delta\binom{i_{s-1}i_s}{1}e^{i_1}(e^0)^{a_1}\cdots e^{i_{s-2}}(e^0)^{a_{s-2}}e^{i_{s-1}}(e^0)^{a_{s-1}+a_s-2n}\\
&+(1-2^{2n})(-1)^{a_s}\binom{2n}{a_s}\delta\binom{i_s}{-1}e^{i_1}(e^0)^{a_1}\cdots (e^0)^{a_{s-2}}e^{i_{s-1}}(e^0)^{a_{s-1}+a_s-2n}                           \\
&+2^{2n}(-1)^{a_s}\binom{2n}{a_s}\delta\binom{i_s}{1}e^{i_1}(e^0)^{a_1}\cdots (e^0)^{a_{s-2}}e^{i_{s-1}}(e^0)^{a_{s-1}+a_s-2n}.
\end{split}
\]
Thus for $(n_1,n_2,\cdots,n_r)\in T_{N,r}$, we have
\[
\begin{split}
&\,\,\,\,\partial_{2n+1}\left(\zeta^{o,\mathfrak{m}}(n_1,\cdots,n_r)  \right)\\
&=\frac{1}{2^r}\sum_{\epsilon_i\in\{\pm 1\},1\leq i \leq r }\epsilon_1\mathcal{I}^{\mathfrak{m}}[\overline{\partial}_{2n+1}\left(e^{\epsilon_1}(e^0)^{n_1-1}\cdots e^{\epsilon_r}(e^0)^{n_r-1} \right) ]\\
&=\frac{1}{2}(2^{2n}-1)\delta\binom{n_1}{2n+1}\zeta^{o,\mathfrak{m}}(n_2,\cdots,n_r)                           \\
&+\frac{(-1)^{n_1-1}}{2}(2^{2n}-1)\binom{2n}{n_1-1}\zeta^{o,\mathfrak{m}}(n_1+n_2-2n-1,n_3,\cdots,n_r)\\
&+\frac{2^{2n}}{2}\delta\binom{n_1}{2n+1}\zeta^{o,\mathfrak{m}}(n_2,\cdots,n_r)                          \\
&-\frac{2^{2n}}{2}(-1)^{n_1-1}\binom{2n}{n_1-1}\zeta^{o,\mathfrak{m}}(n_1+n_2-2n-1,n_3,\cdots, n_r)\\
&+\cdots                          \\
&+\frac{1-2^{2n}}{2}(-1)^{n_{r-1}-1}\binom{2n}{n_{r-1}-1}\zeta^{o,\mathfrak{m}}(n_1,\cdots,n_{r-3},n_{r-2}+n_{r-1}-2n-1,n_r)\\
&-\frac{1-2^{2n}}{2}(-1)^{n_{r-1}-1}\binom{2n}{n_{r-1}-1}\zeta^{o,\mathfrak{m}}(n_1,\cdots,n_{r-2},n_{r-1}+n_r-2n-1)\\
&+\frac{2^{2n}}{2}(-1)^{n_{r-1}-1}\binom{2n}{n_{r-1}-1}\zeta^{o,\mathfrak{m}}(n_1,\cdots,n_{r-3},n_{r-2}+n_{r-1}-2n-1,n_r)\\
\end{split}
\]
\[
\begin{split}
&-\frac{2^{2n}}{2}(-1)^{n_{r-1}-1}\binom{2n}{n_{r-1}-1}\zeta^{o,\mathfrak{m}}(n_1,\cdots, n_{r-2},n_{r-1}+n_r-2n-1)\\
&+\frac{1-2^{2n}}{2}(-1)^{n_r-1}\binom{2n}{n_r-1}\zeta^{o,\mathfrak{m}}(n_1,\cdots,n_{r-2},n_{r-1}+n_r-2n-1) \\
&+\frac{2^{2n}}{2}(-1)^{n_r-1}\binom{2n}{n_r-1}\zeta^{o,\mathfrak{m}}(n_1,\cdots,n_{r-2},n_{r-1}+n_r-2n-1)\\
&=(1-\frac{1}{2})\delta\binom{n_1}{2n+1}\zeta^{o,\mathfrak{m}}(n_2,\cdots,n_r)\\
&-\frac{1}{2}\binom{2n}{n_1-1}\zeta^{o,\mathfrak{m}}(n_1+n_2-2n-1,n_3,\cdots,n_r)\\
&+\cdots \\
&+\frac{1}{2}\binom{2n}{n_{r-1}-1}\zeta^{o,\mathfrak{m}}(n_1,\cdots,n_{r-3},n_{r-2}+n_{r-1}-2n-1,n_r)\\
&-\frac{1}{2}\binom{2n}{n_{r-1}-1}\zeta^{o,\mathfrak{m}}(n_1,\cdots,n_{r-2},n_{r-1}+n_r-2n-1)\\
&+\frac{1}{2}\binom{2n}{n_r-1}\zeta^{o,\mathfrak{m}}(n_1,\cdots, n_{r-2},n_{r-1}+n_r-2n-1)\\
&=(2^{n_{1}-1}-\frac{1}{2})\delta\binom{n_1}{2n+1}\zeta^{o,\mathfrak{m}}(n_2,\cdots,n_r)\\
\end{split}
\]
\[
\begin{split}
&+\frac{1}{2}\sum_{i=1}^{r-1}\left(\binom{2n}{n_{i+1}-1}-\binom{2n}{n_i-1} \right)\\
&\;\;\;\;\;\;\;\;\;\;\;\;\;\; \cdot\zeta^{o,\mathfrak{m}}(n_1,\cdots,n_{i-1},n_i+n_{i+1}-2n-1,n_{i+2},\cdots,n_r ).\\
\end{split}
\]
Thus the proposition holds.
$\hfill\Box$\\

With the help of the above proposition, we can generalize Theorem \ref{5.2} to the case of depth $3$.
\begin{Thm}\label{depth3}
For $r=3$, $N\geq 5\ odd$.\\
(i) The set of the images of elements
$$\{\zeta^{o,\mathfrak{m}}(n_{1},n_{2},n_{3});n_{1}+n_{2}+n_{3}=N,n_{i}\ odd\}$$
in $gr_3^{\mathcal{D}}\mathcal{H}_N$ is a basis of the total space $gr_3^{\mathcal{D}}\mathcal{H}_N$.\\
(ii)Every element in
$$\mathcal{P}^{o}_{N,3}=\langle\zeta^{o}(n_{1},n_{2},n_{3});n_{1}+n_{2}+n_{3}=N, n_{3}>1\rangle_{\mathbb{Q}}$$
can be written as a $\mathbb{Q}$-linear combination of some sum odd multiple zeta values of weight $N$, depth $3$ and multiple zeta values relative to $\mu_2$ of weight $N$, depth less than $3$.
\end{Thm}
\noindent{\bf Proof:}
We have known that the set of elements $$\{\zeta^{o,\mathfrak{m}}(n_{1},n_{2});n_{1}+n_{2}=k,n_{i}\geq 1,\ odd\}$$ is a basis of the space $gr_2^{\mathcal{D}}\mathcal{H}_N$. Similar to the proof of Theorem \ref{5.2} and Theorem \ref{span} we will use the above proposition to prove the first part.  Using Lemma \ref{dy},  we only need to prove that for any given $(n_{1},n_{2},n_{3})\in T_{N,3}$,
\[
\sum_{(k_1,k_2,\cdots,k_r)\neq(n_1,n_2,\cdots,n_r)}|e\binom{k_1,k_2,\cdots,k_r}{n_1,n_2,\cdots,n_r}|<|e\binom{n_1,n_2,\cdots,n_r}{n_1,n_2,\cdots,n_r}|.
\]
 When $n_1=1$, the above inequality is trivial.\\
When $n_{1}\geq 3$, we have
\[
\begin{split}
&\;\;\;\sum_{(k_{1},k_{2},k_{3})\in S_{N,3}}|\delta\binom{k_{3}}{n_{3}}\left({n_{1}-1 \choose k_{1}-1}-{n_{1}-1 \choose k_{2}-1}\right)   \\
&\;\;\;\;\;\;\;\;\;\;\;\;\;\;\;\;\;\;\;\;\;\;\;\;\;\;\;\;\;\;\;\;\;\;\;\;\;\;\;\;\;\;\;\;\;\;\;\;\;\;\;\;\;\;\;\;\;\;\;\;\;\;\;\;\;\;\;+\delta\binom{k_{1}}{n_{2}}\left({n_{1}-1 \choose k_{2}-1}-{n_{1}-1 \choose k_{3}-1}\right)|\\
\end{split}
\]
\[
\begin{split}
&\leq \sum_{(k_{1},k_{2},k_{3})\in S_{N,3}}|\delta\binom{k_{3}}{n_{3}}\left({n_{1}-1 \choose k_{1}-1}-{n_{1}-1 \choose k_{2}-1}\right)|\\
&\;\;\;\;\;\;\;\;\;\;\;\;\;\;\;\;\;\;\;\;\;\;\;\;\;\;\;\;\;\;\;\;\;\;\;\;\;\;\;\;\;\;\;\;\;\;\;\;\;\;\;\;\;\;\;\;\;\;\;\;\;\;\;\;\;\;\;+|\delta\binom{k_{1}}{n_{2}}\left({n_{1}-1 \choose k_{2}-1}-{n_{1}-1 \choose k_{3}-1}\right)|\\
\end{split}
\]
\[
\begin{split}
&\leq \sum_{(k_{1},k_{2},k_{3})\in S_{N,3}}\delta\binom{k_{3}}{n_{3}}\left(|{n_{1}-1 \choose k_{1}-1}|+|{n_{1}-1 \choose k_{2}-1}|\right)\\
&\;\;\;\;\;\;\;\;\;\;\;\;\;\;\;\;\;\;\;\;\;\;\;\;\;\;\;\;\;\;\;\;\;\;\;\;\;\;\;\;\;\;\;\;\;\;\;\;\;\;\;\;\;\;\;\;\;\;\;+\delta\binom{k_{1}}{n_{2}}\left(|{n_{1}-1 \choose k_{2}-1}|+|{n_{1}-1 \choose k_{3}-1}|\right)-2\\
&\leq 4\sum_{i\geq 0}{n_{1}-1 \choose 2i}-2< 2^{n_{1}}-1.
\end{split}
\]
Thus the first statement holds.

As for the second part of this theorem, denote by
$$\mathcal{C}=\{\zeta^{o,\mathfrak{m}}(n_{1},n_{2},n_{3});n_{1}+n_{2}+n_{3}=N\}\setminus\{\zeta^{o,\mathfrak{m}}(n_{1},n_{2},n_{3});n_{1}+n_{2}+n_{3}=N,n_{i}\ odd\}.$$
Assume that there is a lexicographical order on $T_{N,r}$, it induces an order on $\mathcal{C}$ and $\{\zeta^{o,\mathfrak{m}}(n_{1},n_{2},n_{3});n_{1}+n_{2}+n_{3}=N\}$. Let $\alpha$ (resp. $\beta$) be the column vector whose $i^{th}$ element is the $i^{th}$ element in $\mathcal{C}$ (resp. $\{\zeta^{o,\mathfrak{m}}(n_{1},n_{2},n_{3});n_{1}+n_{2}+n_{3}=N\}$). The argument above and Proposition \ref{tran} provide there is a matrix $P$ and an invertible matrix $Q$ such that
$$\partial(\alpha)=P\gamma,\ \partial(\beta)=Q\gamma,$$
where $\gamma=(\zeta^{\mathfrak{m}}(\overline{1})\otimes\zeta^{o,\mathfrak{m}}(1, N-2),\cdots,\zeta^{\mathfrak{m}}(\overline{N-2})\otimes\zeta^{o,\mathfrak{m}}(1, 1))^{T}$.

The last column of  $Q$ is $(0,\cdots,0,2^{N-3}-\frac{1}{2})^T$ obviously, and the last column of $P$ is $0$ because of the following equation:
$$\partial_{N-2}(\zeta^{o,\mathfrak{m}}(n_{1},n_{2},n_{3}))=0,\ \forall \zeta^{o,\mathfrak{m}}(n_{1},n_{2},n_{3})\in \mathcal{C}.$$
By the injectivity of $\partial$ we have
$$\alpha=PQ^{-1}\beta,$$
and that the last raw of $PQ^{-1}$ is $0$. Thus the theorem holds.
$\hfill\Box$\\

Furthermore we can put forward the following conjecture:

\begin{conj}
For any $r\geq 3$, $N\geq r+2$, $N-r\equiv 0\ mod\ 2$, the order $|T_{N,r}|$ matrix
$$E=\left(e\binom{k_1,k_2,\cdots,k_r}{n_1,n_2,\cdots,n_r}\right)$$
as in Proposition \ref{tran} is invertible.
\end{conj}

\begin{rem}
If this conjecture is true we can directly generalize Theorem \ref{depth3} to cases of higher depth by induction.
\end{rem}

\section*{Acknowledgement}
We express our sincere gratitude to Koji Tasaka for pointing out an error in the early version of this manuscript.


\begin{thebibliography}{99}
\bibitem{ref1}F. Brown, \emph{Depth-graded motivic multiple zeta value}, arXiv:1301.3053.
\bibitem{ref2}F. Brown, \emph{Mixed Tate motives over $\mathbb{Z}$}, Ann. of Math, 175(2) (2012), 949-976.
\bibitem{ref3}P. Deligne and A.B.Goncharov, \emph{Groupes fondamentaux motiviques de Tate mixte}, Ann. Sci. École. Normal. Sup. 38(2005), 1-56.
\bibitem{ref4}P. Deligne, \emph{Le groupe fondamental unipotent motivique de $\mathbb{G}_{m}-\mu_{N}$, pour $N=2,3,4,6$ ou $8$}, Publications Mathématiques de l'IHÉS, 112(1) (2010), 101-141.
\bibitem{en}
   B. Enriquez, P. Lochak,
   \emph{Homology of depth-graded motivic Lie algebras and koszulity},
   arXiv: 1407. 4060.
\bibitem{kan}
    H. Gangl, M. Kaneko, D. Zagier,
    \emph{Double zeta values and modular forms},
    Automorphic forms and zeta functions,
    In: Proceedings of the conference in memory of Tsuneo Arakawa,
    World Scientific (2006),
    71-106.
\bibitem{ref5}C. Glanois, \emph{Motivic unipotent fundamental groupoid of  $\mathbb{G}_{m}-\mu_{N}$ for $N=2,3,4,6,8$ and Galois descents}, J. Number Theory, 160 (2016), 334-384.
\bibitem{ref6}J. B. Gil and J. Fres\'{a}n, \emph{Multiple zeta values: from numbers to motives}, Clay Mth. Proceedings, to appear.
\bibitem{hof}
M. E. Hoffman, \emph{The Algebra of Multiple Harmonic Series}, J. of Algebra 194 (1997), 477-495.
\bibitem{ref7}M. Kaneko, K. Tasaka, \emph{Double zeta values, double Eisenstein series, and modular forms of level $2$}, Math. Ann., 357(3) (2013), 1091-1118.
\bibitem{li}J. Li, \emph{Depth-graded motivic Lie algebra},  arXiv:1801.02145v2.
\bibitem{ref8}J. Li, \emph{The depth structure of motivic multiple zeta values},  Math. Ann. (2018). \\https: //doi.org/10.1007/s00208-018-1763-z.
\bibitem{ref9}J. Li, F. Liu, \emph{Motivic double zeta values of odd weight}, arXiv:1710.02244.
\bibitem{ma}D. Ma, \emph{Connections between double zeta values relative to $\mu_N$, Hecke operators $T_N$, and newforms of level $\Gamma_0(N)$ for $N=2,3$}, arXiv:1511.06102.
\bibitem{ref10}I. Soud\`{e}res, \emph{Motivic double shuffle}, Int. J. Number Theory, 6 (2010), 339-370.
\end{thebibliography}
\end{document}